\documentclass[english,a4paper,11pt]{article}
\usepackage{amsmath,amssymb}
\usepackage{graphicx,url}
\usepackage{babel}

\newcommand{\RR}{\mathbb{R}}

\newcommand{\QQ}{\mathbb{Q}}
\newcommand{\ZZ}{\mathbb{Z}}
\newcommand{\CC}{\mathbb{C}}
\newtheorem{Theorem}{Theorem}
\newtheorem{Proposition}[Theorem]{Proposition}

\begin{document}
\title{Elementary fractal geometry. \\ 
3. Complex Pisot factors imply finite type} 
\author{Christoph Bandt\\
Institute of Mathematics, University of Greifswald,  
\url{bandt@uni-greifswald.de}}
\date{\today}
\maketitle

\begin{abstract}
Self-similar sets require a separation condition to admit a nice mathematical structure. The  classical open set condition (OSC) is difficult to verify.  Zerner proved that there is a positive and finite Hausdorff  measure for a weaker separation property which is always fulfilled for crystallographic data. Ngai and Wang gave more specific results for a finite type property (FT), and for algebraic data with a real Pisot expansion factor. We show how the algorithmic FT concept of Bandt and Mesing relates to the property of Ngai and Wang. Merits and limitations of the FT algorithm are discussed.
Our main result says that FT is always true in the complex plane if the similarity mappings are given by a complex Pisot expansion factor $\lambda$ and algebraic integers in the number field generated by $\lambda .$ This extends the previous results and opens the door to huge classes of separated self-similar sets, with large complexity and an appearance of natural textures. Numerous examples are provided. \end{abstract}
 
\section{Overview}\label{intro}
\subsection{Self-similar sets}
A self-similar set is a nonempty compact subset $A$ of $\RR^d$ which is the union of shrinked copies of itself, as defined by Hutchinson's equation 
\begin{equation} A=\bigcup_{f\in F} f(A) \ .  
\label{hut}\end{equation}
Here $F$ denotes a finite set of contractive similarity mappings, called an iterated function system or IFS. A contractive similitude from Euclidean $\RR^d$ to itself fulfils $|f(x)-f(y)|=r_f|x-y|$ where the constant  $r_f<1$ is called the factor of $f.$ We assume that all maps in $F$ have the same factor $r.$ This allows a simple presentation and provides an analogy with crystallographic groups discussed in Section \ref{nb}. For a given IFS  $F,$ there is a unique self-similar set $A$ which is called the attractor of $F.$  See \cite{Bar,BP,Fal} for details. 

With respect to composition of mappings, $F$ generates a semigroup $F^*$ of contractive similitudes. Let $F^n=\{ f_1f_2\cdots f_n\, | f_i\in F\}$ denote the set of compositions of $n$ mappings of $F,$ all with factor $r^n.$ Then $F^*=\bigcup_{n=1}^\infty F^n .$ From \eqref{hut} immediately follows
\begin{equation} A=\bigcup_{f\in F^n} f(A)  \quad \mbox{ for all } n\ge 1\, .
\label{hutn}\end{equation}
Thus $A$ is composed of smaller and smaller pieces, and can be considered as limit set of $F^*.$ Clearly, $A$ is the closure of the set of fixed points of all $f$ in $F^*.$ 

\begin{figure}[h!t] 
\begin{center}
\quad \includegraphics[width=0.4\textwidth]{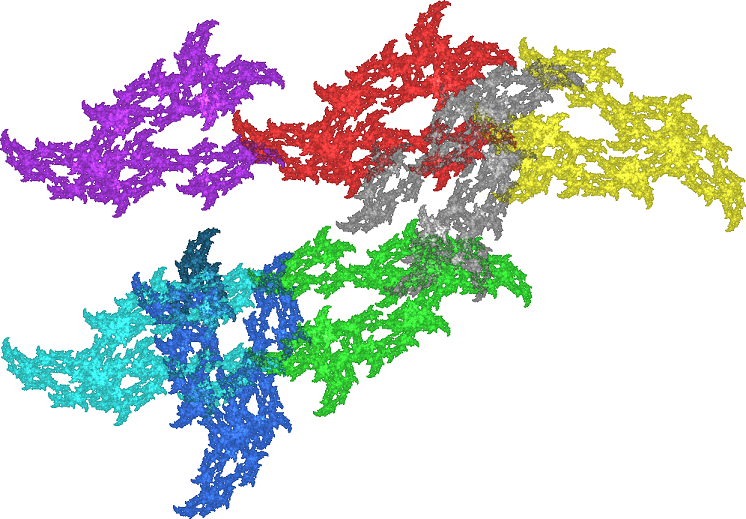} \qquad\quad
\includegraphics[width=0.49\textwidth]{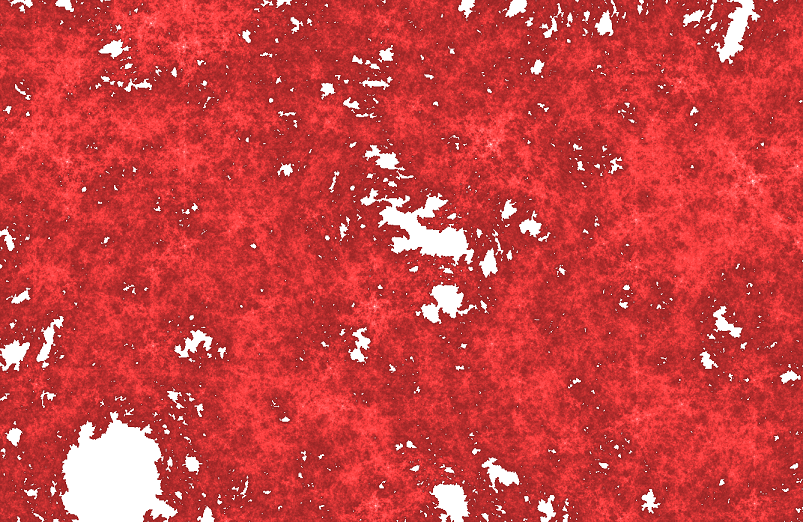} \vspace{.1ex} \\
\includegraphics[width=0.49\textwidth]{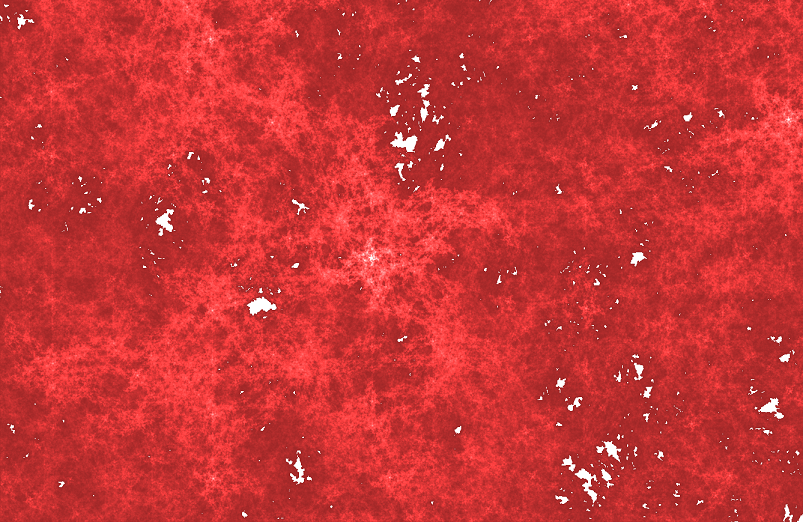} \
\includegraphics[width=0.49\textwidth]{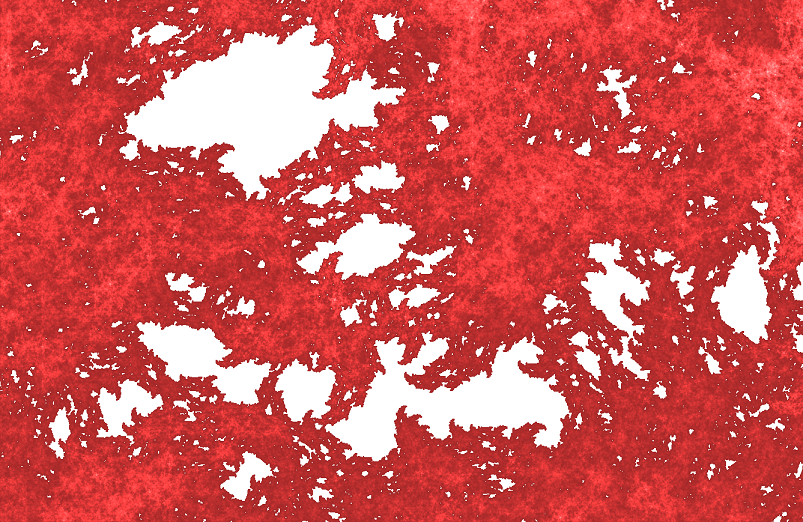} 
\end{center}
\caption{A self-similar set with seven pieces and high complexity constructed from Theorem \ref{main}, with the global view and three particular magnifications. The finite type condition holds with 1799 proper neighbor maps. The data come from the 5th cyclotomic field, as explained in Section \ref{details}. Such sets show an inhomogeneous and random appearance, like natural textures. Due to many overlaps, no holes with the same shape can be seen. Regions with many relatively large holes and with few tiny holes are interlaced. Due to arXiv size restrictions, all figures have been decreased.}\label{Fi1}
\end{figure}  

Self-similar sets appeared as curious examples in topology around 1900, by Cantor, Peano, Hilbert, L\'{e}vy, Sierpinski, and others.
Mandelbrot \cite{FGN} propagated their potential to model natural phenomena, which led Hutchinson to his equation \eqref{hut}.  Barnsley \cite{Bar} developed algorithms to represent real images, like a farn leaf or a forest, as attractors of IFS, in order to compress images.  Meanwhile, self-similar and self-affine sets and measures are well established as theoretical objects in geometrical measure theory, ergodic theory, numeration and related fields \cite{BG,BP,Fal,Fra}.
The hope to model particular sets in nature as attractors of IFS was not fulfilled.  Nevertheless, local views as in Figure \ref{Fi1} may model natural texture in two dimensions - corroded surfaces, images of clouds, soil, tissue, cell cultures etc.  Self-similar sets in this paper have larger complexity and greater affinity to real phenomena than we would expect from known examples.

\subsection{Program of the paper}
Self-similar sets require a separation condition to admit a nice mathematical structure. After a review of strong and weak separation, we introduce the finite type condition (FT) which comes with an algorithm which analyses the structure of $A$ directly from the IFS data. There are various definitions of FT in the literature which will be compared in Section \ref{finit}.

The FT algorithm is given in Sections \ref{fin} and \ref{veri}.  For IFS with algebraic integer data the algorithm is not affected by numerical errors, as proved in Section \ref{nonum}. This makes large IFS and attractors $A$ of high complexity as Figure \ref{Fi1} accessible to computer treatment. The structure of such sets cannot be understood from the global appearance of $A.$  We have to consider many magnifications of $A.$ Examples of similar complexity are known from tilings in Mercat and Akiyama \cite{MA18}.

Our main result says that for large classes of IFS the FT property comes for free. For a complex Pisot expansion factor and algebraic integer data in the IFS, we get FT and weak separation. This is stated in Section \ref{mainresult} and proved in Section \ref{mainproof} by a standard Pisot argument. Section \ref{details}, based on interactive studies with Mekhontsev's IFStile package \cite{M}, presents interesting examples from those classes of IFS.

\subsection{Open set condition} \label{sOSC}
An IFS $F$ fulfils the OSC if there is an open set $U$ such that the $f(U)$ with $f\in F$ are pairwise disjoint subsets of $U.$ The idea is that overlaps $f(A)\cap g(A)$ in this case must be small since they are contained both in the boundary of $f(U)$ and in the boundary of $g(U).$
The OSC was introduced in 1946 by Moran \cite{Mo} to prove that the uniform measure on $A$ is the normalized $\alpha$-dimensional Hausdorff measure. Here $\alpha$ is the similarity dimension of $A$ given by $m r^\alpha =1,$ where $m$ denotes the number of maps in $F.$  The uniform measure assigns to each of the $m$ sets $f(A)$ in \eqref{hut} the value $1/m,$ and to each of the $n$-th level pieces in \eqref{hutn} the value $m^{-n},$ for $n=2,3,...$. With OSC, this natural volume function fits the homogeneous metric structure.

Moreover, we can define clustering numbers $c_n$ for the $n$-th level pieces which have diameter $a_n=r^n\cdot {\rm diam}\, A.$ For $x\in A$ let $c_n(x)$ be the number of  pieces $f(A)$ with $f\in F^n$ which intersect the open ball with center $x$ and radius $a_n.$ Let $c_n=\max \{ c_n(x)|\, x\in A \} .$   In the presence of the OSC, $c=\sup c_n$ is finite so that at some level, the density of pieces does not increase further. Schief \cite{Sch} proved the converse: $c$ is infinite in the absence of OSC. Magnification then provides further clustering of pieces, and the $\alpha$-dimensional Hausdorff measure is zero on $A.$ In summary, OSC is necessary and sufficient to get a nice mathematical structure on $A$ with respect to the similarity dimension $\alpha .$ 

\subsection{Neighbor maps} \label{nb}
An algebraic equivalent for the OSC was given in \cite{BG}, see also \cite{Ba97,Ba00}. A neighbor map for the IFS $F$ has the form $h=f^{-1}g$  with $f,g\in F^*.$ The neighbor set $h(A)$ has the same relative position with respect to $A$ as $g(A)$ has with respect to $f(A),$ up to similarity.  Now $F$ fulfils the OSC if and only if there is a neighborhood $V(id)$ of the identity map in the space of similitudes on $\RR^d$ (with the topology of pointwise convergence) which contains no neighbor maps.

Since we assumed that all mappings in $F$ have factor $r,$ neighbor maps $f^{-1}g$ with $f$ and $g$ in different levels $F^n$ are far from the identity map. It suffices to study the isometries $h=f^{-1}g$ with $f,g\in F^n$ for the same $n.$ In other words, we consider neighbors $h(A)$ which are congruent to $A.$  We call $h$ a proper neighbor map if $h$ is an isometry and $A \cap h(A)\not=\emptyset .$ 

\begin{figure}[h!t] 
\begin{center}
\quad \includegraphics[width=0.26\textwidth]{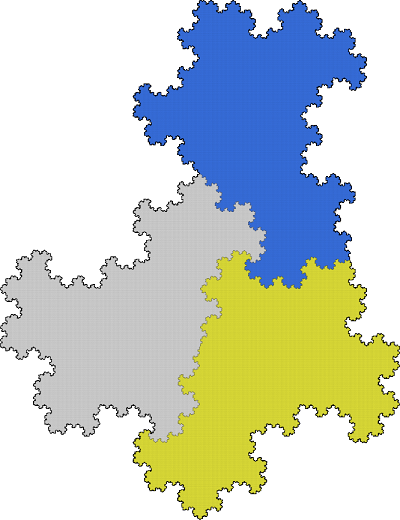} \qquad
\includegraphics[width=0.62\textwidth]{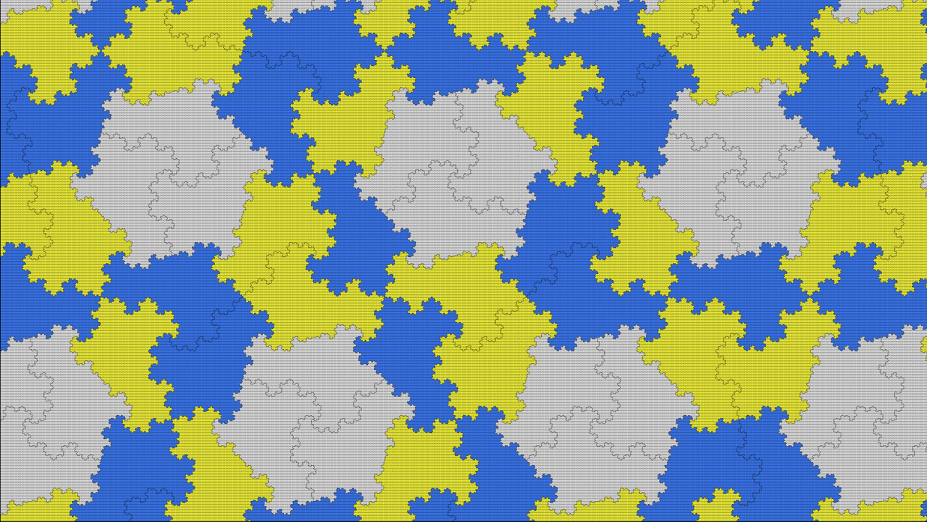} 
\end{center} 
\caption{An aperiodic self-similar tile with three pieces. The generated tiling is colored according to the piece on lowest level. There are 11 proper neighbor maps, including two cases where the intersection is a singleton. Two $120^\circ$ rotations apply to all tiles and generate fractal hexagons. Gray tiles have 5 neighbors, yellow tiles have 7, and blue tiles have 8 neighbors.}\label{Fi2}
\end{figure}  

The classical case is a crystallographic group $H$ with fundamental domain $A.$ In the plane, there are 4, 6, or 8 proper neighbors \cite{GS}. The neighbor maps describe how we can obtain the crystallographic tiling by adding puzzle pieces. This makes sense even without self-similarity. For self-similar crystallographic tilings, see \cite{GG,Lo,LZ17}. Then there are aperiodic tilings as in Figure \ref{Fi2} which have more proper neighbors, but not all possible neighbors apply to each tile.  The generalization to fractal patterns with holes was called fractal tiling by Barnsley and Vince \cite{BV14}, imagining that the construction is extended outwards over the whole plane or space. The local views obtained by zooming-in and by outwards construction are the same, up to similarity. 

\subsection{Weak separation} 
Zerner \cite{Zer} noted that the attractor $A$ can still have a nice mathematical structure when exact overlaps of pieces $f(A)=g(A)$  for $f,g\in F^*$ are permitted in the Bandt-Graf criterion.  The weak separation condition (WSC) is fulfilled if there is a neighborhood $V(id)$ of the identity map in the space of similitudes on $\RR^d$ which contains no neighbor maps except the identity map.  Based on work of Falconer \cite{fal89}, Zerner proved that in the presence of exact overlaps $f^{-1}g=id$ the WSC implies the existence of a positive number $\beta$ smaller than the similarity dimension $\alpha ,$ such that $A$ has Hausdorff dimension $\beta$ and positive finite $\beta$-dimensional Hausdorff measure.

Lau and Ngai \cite{LN99} defined the WSC in slightly different form to analyse the multifractal structure of one-dimensional overlapping self-similar measures, following Lalley \cite{Lal97}.  One-dimensional WSC attractors and measures have been further investigated by Feng \cite{Feng05,Feng16}, Hare and Hare with coauthors \cite{HHM16,HHS18,HR22} and others \cite{dajani21,Kong21}. There are few studies of two-dimensional WSC examples \cite{BM09,BMS,Wu22}, and apparently no studies for greater dimension.

\subsection{Crystallographic data} \label{cyc}
For the case of equal factors, it is convenient to write equation \eqref{hut} like a numeration system. We choose the fixed point of one map  $f_0\in F$ as origin of our coordinate system, define $g=f_0^{-1}$ and $h_k=g\cdot f_k$ for the other elements  $f_1,...,f_{m-1}$ in $F.$ Then   
\begin{equation} g(A)= \bigcup_{k=0}^{m-1} h_k(A)\, ,  \quad  h_k(x)=s_k(x)+v_k
\label{num}\end{equation}
where $g$ is an expanding linear similitude with factor $1/r,$ the $s_k$ are linear isometries, and the $v_k$ are vectors in $\RR^d.$  The isometries $h_k$ are called digits, $h_0(x)=x$ is the identity map and called digit 0. The digits are neighbor maps, $h_k=f_0^{-1}f_k.$  For the decimal system $g(x)=10x, \ A=[0,1],$ and $h_k(x)=x+k, k=0,...,9.$ 

Let $L$ be the lattice generated by the vectors $v_k.$ We say that the IFS $F$ has crystallographic data if $L$ is discrete, isomorphic to $\ZZ^d,$ the mappings $g$ and $s_k$ map $L$ into itself, and the $s_k$ generate a group $S$ which commutes with $g,$ that is, $gS=Sg.$  In this case, the OSC is fulfilled if the $v_k+g(L)$ are disjoint \cite{Ba5,LW1,GG,Lo}.  The expanding map $g$ need not be a similitude. 

\begin{Theorem} (Zerner \cite[Proposition 3]{Zer}) Every IFS with crystallographic data fulfils the WSC. \label{zern}
\end{Theorem}

\subsection{Two examples} 
To clarify the notion of crystallographic data, we give two examples in the complex plane where $g(z)=\lambda z$ and $h_k(z)=s_kz+v_k$  with $\lambda, s_k, v_k\in \CC$ and $|s_k|=1.$ In Figure \ref{Fi2} we consider the hexagonal lattice generated by $1$ and $s=(1+i\sqrt{3})/2.$  This IFS, contained in Mekhontsev's collection of examples \cite{M} since 2016, consists of crystallographic data
\[  g(z)=(s+1)z\, ,\quad h_0=id\, , \quad h_1(z)=-z+1\, , \quad  h_2(z)=-sz -s\, .\]
The OSC is fulfilled, the expanding factor is $\sqrt{3},$ and we have a tile.  All associated self-similar tilings are non-periodic, as in \cite[chapter 11]{GS}.  

 \begin{figure}[h!t] 
\begin{center}
\includegraphics[width=0.49\textwidth]{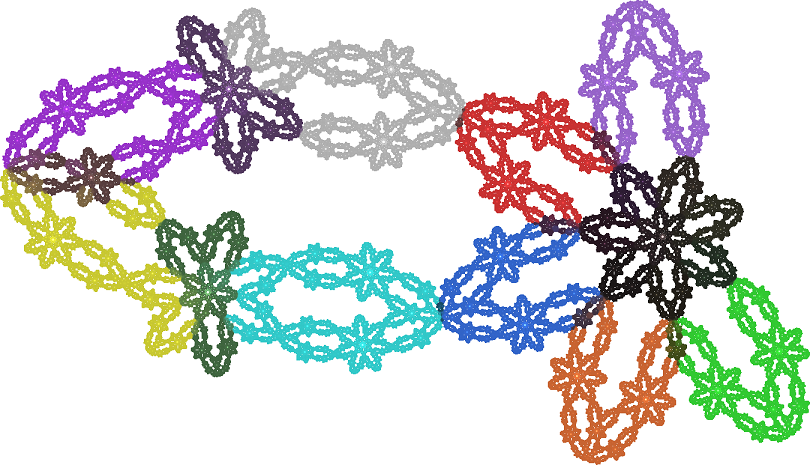} \
\includegraphics[width=0.49\textwidth]{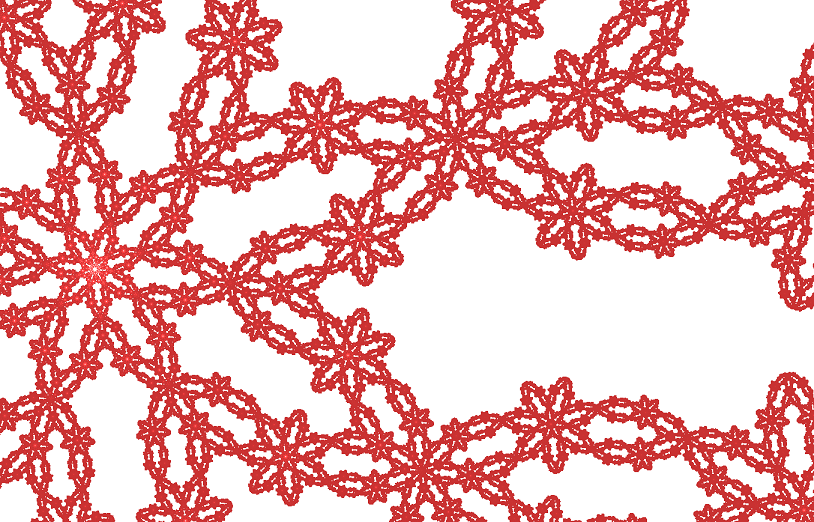} 
\end{center}
\caption{An overlapping attractor with nine pieces and sevenfold symmetry. All 13 proper neighbor maps contain exact overlaps and can be identified in the figure. They form two sets of six rotations around the center of the rosette and around the left tip of the set, plus the identity map for complete overlap. In the local view on the right, we can walk from each piece to any other by applying rotations between neighboring pieces, either around their rosette center or around their tips, which become centers of incomplete rosettes. }\label{Fi3}
\end{figure}  

An attractor with $m=9$ strongly overlapping pieces is shown in Figure \ref{Fi3}. The data are non-crystallographic since they contain rotations around $\beta=\pi/7.$  The number $s=\cos \beta+ i\sin\beta$ fulfils $s^7=-1.$
We take the expansion factor $\lambda= 1+s+s^2\approx 2.5245-1.2157i.$ The origin of our coordinate system is taken as center of the rosette on the right of the global view. The mappings $h_1,...,h_5$ are $z, s^2z, s^{-2}z, s^4z, s^{-4}z.$ The corresponding pieces form the rosette and create exact overlaps of second-level pieces. They generate six rotations $s^{2k}, k=1,...,6$ around zero as neighbor maps.

A second rotation center for neighbor maps is the tip of $A$ on the left, which becomes the center of incomplete rosettes in the local view on the right.
The parallel pieces in the middle of $A$ are defined by $h_6=s^8(z+s^3)$ and  $h_7=s^8(z+s^2).$  The two pieces on the left given by $h_8=s^2(z-s)$ and $h_9=z-s^4$ intersect at their tips which coincides with the tip of $A.$ The parallel pieces $A_6$ and $A_7$ also intersect other pieces at their tips. The corresponding neighbor maps are six rotations which act on the incomplete rosettes of the local view. Together with the identity map, which indicates complete overlap of pieces, this gives 13 neighbor maps between intersecting pieces. Applying the two types of neighbor rotations, we can walk from each piece to any other piece of a fractal tiling, as sketched on the right of Figure \ref{Fi3}.

\subsection{The finite type algorithm} \label{fin}
For IFS which do not consist of crystallographic data, the weak separation property is hard to check. We use the following property \cite{BM09} which is stronger as proved by Nguyen \cite{nguyen}. \vspace{2ex}

\noindent {\bf Definition. }  An IFS $F$ is of finite type (FT) if there are finitely many proper neighbor maps.\vspace{2ex}

The first, slightly different and  more technical concept of FT was introduced by Ngai and Wang \cite{NW}, extending work of Lalley \cite{Lal97}. Our definition says that there are finitely many pairs of intersecting pieces, up to similarity. In essence, the definition of Ngai and Wang says that there are finitely many pairs of overlapping pieces, in the sense that some subpieces of the two pieces coincide, up to similarity. Thus our property is a little bit stronger.  We discuss the two concepts carefully in Section \ref{finit}. 

Our finite type condition can be checked by an algorithm sketched in \cite{Ba97,Lal97,Zer}, more carefully described in \cite{Ba00,LNR}, and for lattice tilings in \cite{ScT}. Here we continue the discussion of this method. The algorithm seems the only way to analyze complicated examples from given IFS data, even when we are only interested in overlapping pieces. In Section \ref{nonum} we show that for IFS with algebraic integer data, the FT algorithm is not affected by numerical errors.  It works accurately in the presence of thousands of neighbor maps. 

When the algorithm verifies the FT, it also decides the OSC and provides a finite automaton which we called neighbor graph. In the case of OSC, this automaton completely describes the topology of the attractor and of its dynamical boundary, as well as the dimensions of the boundary \cite{EFG1,EFG2,BM09}.  The above definition and the algorithm extend in a straightforward way to the case of different contraction factors as well as to graph-directed systems. All this is implemented in Mekhontsev's IFStile package \cite{M}. \vspace{2ex}

{\bf The FT algorithm:  Deciding FT.  } On the one hand, proper neighbor maps  have the form 
\begin{equation}  h=f_{n}^{-1}\cdots f_{1}^{-1}g_{1}\cdots g_{n} \ \mbox{ with } \ f_{k},g_{k}\in F \ \mbox{ for } k=1,...,n \, . \label{nei} 
\end{equation} 
On the other hand, 
$h(x)=s(x)+v$ with a linear isometry $s$ and a vector $v=v(h).$ There is a constant $C,$ calculated in Section \ref{veri} from the data of $F,$ such that $h$ is not proper if $|v(h)|>C.$ The proper neighbor maps are now generated by recursion over $n$ in \eqref{nei}. We start with $H_0=\{ id\}$ and $H_1= \{ h_1=f^{-1}g\, |\, f\not= g\in F, |v(h_1)|\le C\}.$

Suppose we have constructed $H_k$ for $1\le k\le n.$  Then let
\[ H_{n+1}=\{ h_{n+1}=f^{-1}h_ng \mbox{ with } h_n\in H_n, \ f,g\in F\mbox{ and }|v(h_{n+1})|\le C\} \setminus \bigcup_{k=0}^n H_k\ .\] If $H_{n+1}$ is empty, FT is fulfilled. Otherwise we go to the next level.

This method must be implemented in exact arithmetics, which can be done for algebraic data, as explained in Section \ref{nonum}. In the presence of numerical errors, we would not be sure whether $h_{n+1}$ coincides with some $h_m$ of a previous level.  If FT holds, the algorithm will confirm FT in finite time, see Section \ref{veri}. The algorithm cannot safely reject FT in finite time. Nevertheless, it works surprisingly well in finding simple examples with few neighbor maps, say less than 100. In that case we stop the algorithm when we have reached a much larger number of neighbor candidates, say 20000, concluding that it is unlikely that there are less than 100 neighbor maps. Then we try another randomly chosen IFS \cite{M}. \vspace{2ex}

{\bf The FT algorithm:  Finding proper neighbors.  }
For this purpose we set a link from each $h_n$ to each possible successor $h_{n+1}=f^{-1}h_ng.$  When FT is confirmed at level $n^*+1,$ we erase all $h_n$ without successor in $\bigcup_{k=0}^{n} H_k$ and their incoming links, starting with $n=n^*$ and going down up to $n=1.$  The proper neighbor maps which remain will admit an outgoing path of links to a cycle of links. These links form the neighbor graph. See \cite{Ba00,EFG1,EFG2,BM09,HR22,Lo,LZ17,ScT,TZ20} for details and examples. Overlap neighbors are those which admit an outgoing path to the initial vertex $id.$ If they do not exist, the OSC holds.

\subsection{Fractal constructions and Pisot numbers} \label{fracpis}
Algebraic numbers come into IFS since an exact overlap $f=f'$ for some $f,f'\in F^n,$ and a cycle of links in the above algorithm, as well as other geometric properties of $A,$ are defined by polynomial equations. An algebraic integer $\lambda$ is the root of a polynomial 
\begin{equation}
P(\lambda)=\lambda^d+\sum_{k=0}^{d-1} m_k\lambda^k\qquad \mbox{ with integer coefficients } m_k.
\label{Poflambda}\end{equation}
We take the smallest possible $d$ which is called the degree of $\lambda .$ A real root $\lambda >1$ is a Pisot number if all other roots of $P$ have modulus smaller than one. A complex root with $|\lambda |>1$ is called complex Pisot number if all roots of $P$ except $\lambda$ and $\overline{\lambda}$ have modulus smaller than one. 

We note some key results on fractals and Pisot numbers. The simplest one-dimensional IFS are Bernoulli convolutions, with $g(z)=\lambda z$ and $h_0(z)=z, h_1(z)=z+1.$ In 1939, Erd\"os \cite{E39} discovered the special r\^ole of Pisot numbers $\lambda$ for the Fourier transform of the self-similar measure. Only in 1962, this was clarified by Garsia \cite{Ga} and led to a lot of research in one-dimensional self-similar measures, cf. \cite{AFH,Feng05,Feng16,FS92,HHM16,HR22,Lal97,LN99}.  In 1982, Rauzy found a connection between substitutions with a matrix with Pisot eigenvalue, dynamical systems and self-similar tilings \cite{Rauzy} which triggered much work on fractal tilings \cite{A02,AI,MS14,RWY14} and the still unsolved Pisot substitution conjecture \cite{ABB,BK06,BS05,MA18}. 
The Penrose self-similar tilings \cite{GS}, which became models of quasicrystals, are also generated by Pisot numbers. Research in aperiodic tilings \cite{Baake2013,FG18} just in 2023 led to the discovery of a single aperiodic tile \cite{hat}. In his visionary lecture notes from 1989, Thurston \cite{T89} discussed number systems with Pisot base, corresponding finite-state automata and self-similar tilings generated by complex Pisot numbers. Solomyak \cite{SoDyn} proved in 1997 that the dynamical system generated by a self-similar tiling in one or two dimensions has pure discrete spectrum if and only if the expansion factor is real or complex Pisot and a certain density condition holds which can be decided by an automaton similar to our neighbor graph.

Ngai and Wang \cite{NW} used the finite type property to define a neighborhood graph from which the dimension $\beta$ of an overlapping attractor can be explicitly determined as $\log\sigma /\log \lambda$ where $\sigma$ is the spectral radius of the incidence matrix $M$ of the graph, and $\lambda$ the real expansion constant of the IFS.  They also presented a class of finite type IFS which goes beyond the discrete data in Zerner's theorem. 
For an algebraic integer $\lambda$ of degree $d,$ let
\begin{equation}
\ZZ(\lambda)=\{ z=\sum_{k=0}^{d-1} n_k\lambda^k\ | \ n_k\in\ZZ \} .
\label{algint}\end{equation}

\begin{Theorem} (Ngai and Wang \cite[Theorem 2.5]{NW}) Consider an IFS on $\RR^d$ in the form \eqref{num} with $g(x)=\lambda x$ where $\lambda >1$ is a Pisot number. Assume that the orthogonal matrices $s_k$ generate a finite group $G,$ and that there are real numbers $r_1,...,r_d$ such that for each $g\in G$ and each $v_k$
\[ g(v_k)\in r_1\ZZ (\lambda )\times r_2\ZZ (\lambda )\times ...\times r_d\ZZ (\lambda ) .\]
Then the finite type property is fulfilled. \label{NWtheorem}
\end{Theorem}

This certainly applies to $d=1$ which includes the examples in \cite{NW} and previous results \cite{Ga,Lal97}. We do not discuss $d\ge 3$ where there are few finite isometry groups and so far no examples.  For $d=2,$ it applies to the golden gaskets \cite{BMS} and one gasket example in \cite{NW}, and to rotation groups of order 3, 4, and 6. It is not clear how to choose $r_k$ for rotation of order 5 or 7. The theorem does not apply to Figure \ref{Fi2} and not to the boundary IFS of the Levy dragon (as claimed in \cite{NW}) because $\sqrt{3}$ and $\sqrt{2}$ are not Pisot numbers. Since the complex factors $s+1$ and $1+i$ are complex Pisot numbers, those will be covered by Theorem \ref{mainresult}. The dense product lattice in Theorem \ref{NWtheorem} does not allow expansion maps which involve an irrational angle.

Both Zerner and Ngai and Wang are a bit more general than quoted here, allowing IFS maps with factors  $r^n$ for different $n.$  This could also be done for Theorem \ref{mainresult} below, with more technical effort.  We prefer equal factors for which neighbor maps are isometries.  Moreover, our IFS will not contain reflections which lead to fractals with different appearance. We focus on IFS of conformal maps. Reflections can be included using complex conjugation \cite{M}. 

\subsection{Main result} \label{mainresult}
In the complex plane, we take an IFS  $F=\{ f_k=g^{-1}h_k\, |\, k=1,...,m\}$ in the form 
\begin{equation}
g(A)= \bigcup_{k=0}^{m-1} h_k(A)\, ,  \quad  g(z)=\lambda z\, \mbox{ and } h_k(z)=s_kz+v_k\, .
\label{alg}\end{equation}

\begin{Theorem} \label{main} If $\lambda$ is a complex Pisot number, $s_k, v_k\in \ZZ(\lambda) ,$ and the rotations $s_k$ are rational, the IFS is of finite type.
\end{Theorem}

This simple and natural statement contains all two-dimensional examples of
Theorem \ref{zern}, as well as applications of Theorem \ref{NWtheorem} when we replace a real Pisot number $\lambda$ by the complex factor $s\lambda$ and change the IFS without change of the attractor.  We allow non-crystallographic rotations, but the $s_k$ must satisfy $s^n=1$ for some $n.$ Note that  $\ZZ(\lambda)$ is a dense lattice in $\CC .$ 
The proof in Section \ref{mainproof} uses the Pisot property in a standard way. 

The new feature is our computational approach with the FT algorithm which allows further analysis of the attractors.  Section \ref{details} provides various interesting and highly complex examples for which neighbor graphs were calculated and connectedness properties verified. In Section \ref{newalgo} we sketch how the dimension formula of Ngai and Wang with spectral radius of a matrix can be calculated directly from the IFS data even for thousands of neighborhoods. This is work in progress. 

\subsection{Acknowledgement} 
All examples were found and all figures produced by the wonderful free package of Dmitry Mekhontsev \cite{M}. When we met in autumn 2016, he had completed the first version of IFStile in his spare time. I invited Dmitry to transform part of his work into the thesis \cite{Mth} which he defended in January 2019. He taught me more than I could teach him. But he does not like writing papers, and our cooperation has faded. Dmitry still maintains the IFStile finder which was made for OSC examples and contains algebraic number fields in linear algebra disguise. When I asked him to include exact overlaps, he kindly added a button OVL to his search menu.  I sincerely hope that \cite{EFG1,EFG2} and the present paper will encourage colleagues to use \cite{M}.

\section{The finite type property} \label{finit}
\subsection{Two neighbor concepts}
We compare our concept from \cite {BM09} with the finite type concept of
Ngai and Wang \cite{NW}.  For details on various similar conditions in the literature we refer to the surveys \cite{DE11,DLN13,HHR21}.
In both cases, equivalence classes of certain configurations with respect to similarity are counted. Ngai and Wang considered neighborhood types, Bandt and Mesing studied neighbors $h(A)$ of the standard piece $A.$ This does not matter when only the question of finiteness is considered.

\begin{Proposition}\label{P1}
For any neighbor concept, the number of neighbors is finite if and only the number of neighborhoods is finite.
\end{Proposition}

{\it Proof. } If there are $M$ neighbors, there can be at most $2^M$ neighborhoods. If there are $N$ neighborhoods with $k_1,...,k_N$ neighbors, respectively, there are at most $\sum_{j=1}^N k_j$ neighbors.
\hfill $\Box$ \vspace{1ex}
 
Now let us consider the neighbor concepts. Our definition can be formulated as
\begin{equation} f,g\in F^*\mbox{ are proper neighbors if } f(A)\cap g(A) \not=\emptyset
\label{prop}\end{equation}
The neighbor concept of Ngai and Wang \cite{NW} depends on an open set $V$ which fulfils $f(V)\subset V$ for all $f\in F.$ This implies that $A$ is a subset of the closure of $V.$ The sets $f(V)$ need of course not be disjoint as in the OSC since exact overlaps are allowed. 
\begin{equation} f,g\in F^*\mbox{ are neighbors with respect to }V\mbox{ if } f(V)\cap g(V) \not=\emptyset
\label{neiV}\end{equation}
In general, proper neighbors cannot be realized by a choice of $V.$ However:

\begin{Proposition}\label{P2}
Two maps $f,g\in F^*$ are proper neighbors if and only if they are neighbors with respect to every $V\supset A.$
\end{Proposition}

{\it Proof. } If $f(A)$ and $g(A)$ intersect, this holds also for all their neighborhoods. If the sets do not intersect, there is an $\varepsilon >0$ with $d(f(A),g(A))=\varepsilon$ where $d(X,Y)=\inf\{ |x-y|\,|\, x\in X, y\in Y\}\ .$   Let $V=\{ x|\, d(x,A)<\varepsilon/2 \}.$ Since $f,g$ are contractions, $f(V)\cap g(V)=\emptyset .$   \hfill $\Box$ \vspace{1ex}

The dependence on $V$ is a strength and a weakness of the Ngai-Wang concept. The set $V$ can be adapted to the purpose. Taking a larger set $V,$ we get more neighbors. When drawing a picture of a fractal tiling around $A,$ we should include all neighbors in our window. When calculating the dimension of an overlapping construction, Ngai and Wang took small sets $V$ in order to minimize the size of the corresponding matrix \cite{NW}. For a definition, however, the dependence on an open set $V$ is an obstacle. It is not easy to decide whether two pieces are neighbors, unless $V$ is specifically chosen for the IFS.  

\subsection{Verification of the FT algorithm} \label{veri}
The big advantage of our neighbor concept is that in the FT case all proper neighbor maps can be computed by the FT algorithm in Section \ref{fin} from the IFS data. Let us now calculate the constant $C$ and verify that the algorithm stops in the case that FT is fulfilled.

We write $F$ in the form \eqref{num} with digits $h_k(x)=s_k(x)+v_k.$ Then the closed ball $B$ around the origin with radius  $R= r\cdot \max_k |v_k| / (1-r)$ is mapped by the IFS mappings $f_k=g^{-1}h_k$ into itself, since $r|v_k|\le (1-r)R.$ Thus $A$ is contained in the ball $B.$ If the neighbor map $h$ fulfils $|v|>2R$ then it maps $B$ into a disjoint ball and cannot be proper. This determines the constant $C$ for the FT algorithm and the necessary condition for proper neighbor maps \cite{Lal97,Ba00,NW,LNR}:
\begin{equation}
|v|\le C=2R = \frac{2r}{1-r}\cdot \max_k |v_k| \ .
\label{Ccon}\end{equation}
Now assume that FT holds, so that we have $N$ proper neighbor maps. Then our algorithm stops with $H_{n}=\emptyset$ for some $n\le N$ because the set of neighbor maps increases at each level. Thus our algorithm will stop after finite time. It would be nice to have a quantitative estimate.

The WSC* condition of Lau, Ngai and Rao \cite{LNR} says that \eqref{Ccon} is fulfilled for finitely many neighbor maps. So this is just a technical equivalent of our FT property. 

\subsection{Equivalence of the finite type concepts}
The equivalence of Ngai-Wang FT and WSC* for a special choice of $V$ was shown in \cite[Theorem 6.4]{DLN13}. We now prove that no choice is necessary.

\begin{Theorem} \label{equi}
Let $F$ be an IFS with equal factors on $\RR^d.$
\begin{enumerate}
\item[(i)] If our finite type condition holds, then for every $D>0$ there is a finite number of neighbor maps $h(x)=s(x)+v$ with $|v|\le D.$
\item[(ii)] For any open set $V\supset A ,$ our concept of finite type coincides with the finite type condition of Ngai and Wang with respect to $V.$
\end{enumerate}
\end{Theorem}

{\it Proof. } (i). Assume FT holds and $H_{n}=\emptyset .$ We consider the neighbor maps $h$ in $H=\bigcup_{k=0}^n H_k$ which must fulfil \eqref{Ccon}, plus their immediate successors $h'$ which are outside $H$ and do not fulfil  \eqref{Ccon}. Clearly this is a finite set.  From this set we take the subset $H'$ of maps which are not proper, and define the positive number 
\[ \delta =\min\{ d(A,h(A))\, |\, h\in H' \} \ . \]
Since the algorithm stopped, $H'$ contains all immediate successors of proper maps which are not themselves proper. So each improper neighbor map must be either in $H'$ or a successor of an element of $H'$ over several generations. We now choose an integer $k$ with  $r^{-k}\delta >D.$ Assuming that the origin is in $A$ (which will hold if $h_0=id$), we show that the set of all neighbor maps with $|v|\le D$ is contained in the union $H^*=H\cup H'\cup H^1\cup\dots\cup H^{k-1}$ where $H^j$ denotes the set of all successors of $H'$ in the $j$-th generation. 

The following argument is needed here. When $h=f^{-1}g$ is a neighbor map  and  $h'=f_k^{-1}hf_j$ an immediate successor of $h,$ then $d(A,h'(A))\ge r^{-1}d(A,h(A)) .$
For $f\in F^n$ this follows from 
\[  r^{n+1}\cdot d(A,h'(A))=d(ff_k(A),gf_j(A)) \ge d(f(A),g(A))=r^n\cdot d(A,h(A))\, .\]
Now if $h''$ is not in $H^*,$ it is a successor of $H'$ in generation $k$ or larger. So 
\[ |v(h'')|\ge d(0,h''(A))\ge d(A,h''(A))\ge r^{-k}\delta >D\ .\]
This shows that all neighbor maps with $|v|\le D$ must be in $H^*.$  Since $H^*$ is finite, this implies (i).

For (ii), let $V\supset A.$ The Ngai-Wang neighbors contain the proper neighbors by Proposition \ref{P2}. So the Ngai-Wang FT condition implies our FT condition.  On the other hand, the neighbor maps with $|v|\le D$ contain the Ngai-Wang neighbors whenever $D>2\cdot {\rm diam\,} V\, .$ By (i), our FT implies the FT condition of Ngai and Wang. 
\hfill $\Box$

\subsection{Examples where neighbor concepts differ}
For sets $V$ which do not contain $A,$ some proper neighbors may fail the neighbor property in the Ngai-Wang setting. A simple example is the golden mean Bernoullli convolution: $f(x)=tx, g(x)=tx+1-t$ with $t=(\sqrt{5}-1)/2$ and $A=[0,1].$  On any level $n,$ two pieces will either be disjoint, have one point in common, have a small overlap like $f(A),g(A),$ have a large overlap like $fg(A),gf(A),$ or coincide altogether. There are seven proper neighbor maps: $h(x)=x\pm 1, h(x)=x\pm t, h(x)=x\pm t^2,$ and $h=id.$

It is very natural to choose $V=(0,1).$ Then we get five Ngai-Wang neighbors. The two neighbors with one-point intersection are neglected.  For one-dimensional attractors, this makes sense and is taken as FT definition in recent papers by Hare and coauthors \cite{HHS18,HR22}.

If we take the interval as attractor of two maps with non-commensurable factors, say $f(x)=0.7x$ and $g(x)=0.3x+0.7,$ we will have infinitely many neighbors with one common point, because of the different ratios of the neighbors' lengths. For $V=(0,1)$ there are no Ngai-Wang neighbors. In this case the two FT concepts are different.

\begin{figure}[h!t] 
\begin{center}
\quad \includegraphics[width=0.36\textwidth]{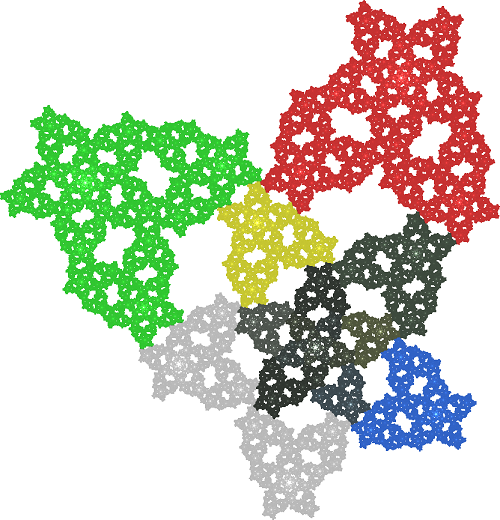} \qquad
\includegraphics[width=0.56\textwidth]{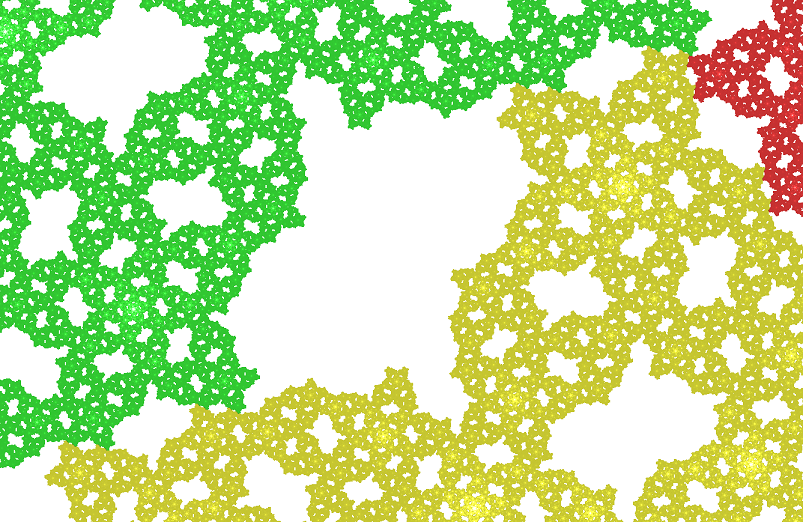} 
\end{center}
\caption{A self-similar set with five pieces, where three pieces meet in the little rosette in the lower part, and two of them also have the black piece of second level in common.  Let $V$ be the interior of the filled-in set. There are two types of Ngai-Wang neighborhoods: the rosette which is the union of three overlapping pieces, and a piece without overlaps. Ngai and Wang interpret the upper two pieces as isolated, and disregard the intersection of pieces in a Cantor set, and the triple intersection in a point which are apparent in the magnification.}\label{Fi4} 
\end{figure}  

A more intricate example is shown in Figure \ref{Fi4}.
The pieces $A_1$ (grey), $A_2$ (yellow) and $A_3$(blue) of the self-similar set $A$ overlap and together form a rosette $R_0.$ The other two pieces $A_4$ (green) and $A_5$ (red) intersect each other and the rosette in Cantor sets. In the Ngai-Wang approach, the interior of the filled-in set can be conveniently taken as open set $V.$ Then the three overlapping pieces together form a neighborhood, while there are four possible pairs of neighbors. The other two pieces are isolated and represent the trivial neighborhood because the filled-in interiors do not meet, for example $f_4(V)\cap f_5(V)=\emptyset .$

It is now possible to construct a big rosette $R$ by adding three more copies $f_6(A), f_7(A), f_8(A)$ cyclically around the rosette $R_0.$ Moreover, $R_0=f_0(R)$ will be a small copy of $R.$ This establishes a so-called graph-directed IFS from the two neighborhood types $R$ and $A$:
\[ R=f_0(R)\cup\bigcup_{k=4}^8 f_k(A) \quad , \quad  A=f_0(R)\cup f_4(A)\cup f_5(A)\ .\]
The second equation describes what you see on the left of Figure \ref{Fi4}. There is no overlap between the pieces, the OSC holds and the Hausdorff dimension can be determined from a matrix equation \cite{NW}.  In this example there is a much simpler way to eliminate the overlaps: $A$ is a non-overlapping ordinary self-similar set with three pieces $A_k$ and three smaller second-level pieces:
\[  A= f_4(A)\cup f_5(A)\cup f_1(A)\cup f_2f_4(A)\cup f_2f_5(A)\cup f_3f_5(A)\ .\]
This directly gives the Hausdorff dimension $\alpha$ by the equation for the similarity dimension $3r^\alpha +3r^{2\alpha} =1$ where $r$ is the factor of the $f_k.$ Setting $x=r^\alpha ,$ the quadratic equation $x^2+x=\frac13 $ yields $x=\frac16 (\sqrt{21}-3)$ and $\alpha =\log x/\log r\approx 1.849\, .$

The data of this IFS were found with the method of Section \ref{prescribed}, and the expanding factor $\lambda\approx 2.0420+0.2384i$ fulfils the equation $\lambda^2=-(s^2+s^3)\lambda -s^3$ where $s$ is the standard 5th root of unity, $s^5=1.$ Thus $\lambda$ defines a quadratic extension of the 5th cyclotomic field and has degree 8. This factor, as well as that of Figure \ref{Fi3}, is not Pisot. In both examples there is only one type of hole. In \cite{EFG2} we constructed a ``dog carpet''. Here we have a ``moth carpet'' and in Figure \ref{Fi3} a ``cello carpet''. One-hole self-similar sets are related to tiles. Both a vertex in a tiling and a hole in a fractal tiling represent a generating relation $h_1h_2...h_n=id$ of neighbor maps.  

While Ngai-Wang neighbors determine the Hausdorff dimension $\beta$ of $A,$ the proper neighbors determine details of the structure: find out which pairs of pieces have points in common (all except (1,5) and (3,4)), what are the Hausdorff dimensions $\gamma$ of their intersections (all have the same $\gamma=\log 2/(-2\log r) \approx 0.481$) and so on.

\subsection{Restricted neighbor concepts}
It can make sense to neglect certain proper neighbors.  Another example is any crystallographic tiling of the plane by a triangle $A,$ generated by reflections \cite{GS}. There are three neighbors $h(A)$ which share an edge with $A.$ The corresponding reflections $h$ generate the Coxeter group. Numerous neighbors which share only one point with $A$ are not so interesting. Ngai and Wang would exclude them by taking an open set $V$ such that $A\setminus V$ consists of the three vertices of $A.$ We prefer to say that we study edge neighbors and neglect the point neighbors \cite{BM09}.

A general way to define special neighbors is to take a set function $\Phi$ and a number $c$ and require that \[ \Phi(f(A)\cap g(A)) \ge c\ .\]
Here $\Phi$ can be a cardinality, a dimension, or a measure.   For Ngai and Wang, $\Phi$ is the Hausdorff dimension and $c={\rm dim}\, A.$ Under the assumption of finite type, this means that $f(A)\cap g(A)$ contains an exact overlap, that is, $ff'=gg'$ for some $f',g'\in F^*.$ The most restrictive neighbor concept is $f(A)=g(A),$ studied by Neunh{\"a}userer \cite{neun}.  
We now show how the FT algorithm allows to abandon the sets $V.$

\subsection{From open sets to computer algorithms}\label{newalgo}
As noted in Section \ref{fracpis}, Ngai-Wang neighbors were a tool to explicitly determine the dimension of WSC attractors. The formula with the spectral radius of the neighborhood incidence matrix $M$ greatly improved Zerner's approximation method.  In dimension one, $M$ is easily determined when the set $V$ is an interval \cite{HR22}. In dimension greater one, however, the graph of neighborhoods of an open set $V$ is very difficult to establish, unless $A$ is as simple as Figure \ref{Fi4}. We suggest to determine $M$ in the following algorithmic way. 

\begin{enumerate}
\item Determine all proper neighbors and establish the neighbor graph,
\item Reduce to those neighbors which include overlaps (the vertices with a path to $id$),
\item Determine the neighborhood graph of the reduced neighbor graph,
\item Calculate the incidence matrix $M$ and its spectral radius.
\end{enumerate}

All these steps can be formalized and delegated to the computer. This is work in progress.
The method requires that the IFS fulfils FT in our sense.  There is an example by Kenyon with  infinitely many proper neighbor maps and OSC and therefore without overlaps  \cite[Example 3.8]{DLN13}. For this construction, the method fails. However, if an IFS is just given by data, and we look for overlaps without geometric information, we have to check all proper neighbors, since any small intersection $A\cap h(A)$ could contain a piece of $A.$ In this sense, the two FT concepts again coincide.

\section{IFS in algebraic number fields} \label{proofs}
\subsection{Basic setting}\label{sett}
We now explain that the FT algorithm for algebraic data works with integer vectors only, and then prove Theorem \ref{main}.  The expanding map of the IFS \eqref{alg} is  $g(z)=\lambda z.$   The number $\lambda$ should be a non-real algebraic integer with modulus larger than one, but need not be a Pisot number yet. Thus $\lambda$ is a root of the minimal polynomial $P$ in \eqref{Poflambda} of degree $d.$ The lattice $\ZZ(\lambda)$ of algebraic integers is given in \eqref{algint}. Replacing integer coefficients $n_k$ by rational numbers $x_k,$ we obtain 
\begin{equation}
\QQ(\lambda)=\{ z=\sum_{k=0}^{d-1} x_k\lambda^k\ | \ x_k\in\QQ \} \subset\CC .
\label{algrat}\end{equation}
The algebraic number field $K$ generated by $P$ is in principle this set, only considered as $d$-dimensional vector space over $\QQ$ with basis $B=\{ 1,\lambda ,\lambda^2,..., \lambda^{d-1} \} .$ The complex number $z\in\QQ(\lambda)$ is replaced by the column vector $\tilde{z}=x=(x_1,...,x_d)'$ in $K.$ Cf. \cite{cohen,Baake2013}. Let $\overline{K}\supset K$ denote the corresponding real vector space with coefficients $x_k\in\RR .$ An element of $\ZZ(\lambda)$ corresponds to an integer vector in $K$ and $\overline{K}.$

Multiplication with an element of $K$ is a linear map on $K,$ given by a $d\times d$ matrix with respect to the basis $B.$  Multiplication by $\lambda$ is given by the companion matrix $L$ of $P,$ with first two columns $(0,1,0...0)',(0,0,1,0...0)'$ because $1\cdot\lambda=\lambda , \ \lambda\cdot\lambda=\lambda^2,$ and last column $(-a_0,-a_1,...,-a_{d-1})'$ because $P(\lambda)=0$ in \eqref{algint}.  

The digits of the IFS \eqref{alg} are $h_k(z)=s_k z + v_k, k=1,...,m$ where $s_k$ and $v_k$ are in $\ZZ(\lambda)$ and $s_k^n=1$ for some $n.$ As additive elements $s_k,v_k$ correspond to $d$-dimensional integer vectors $\tilde{s_k},\tilde{v_k}$ in $K$ which come from their representation \ref{algint} with respect to the basis $B.$ For specific examples, it is not easy to find this representation, and sometimes another basis may be more convenient, as discussed in Section \ref{details}. 

Anyway, when we can express $s_k$ in terms of the basis, this holds also for $s_k\lambda , s_k\lambda^2$ etc., and we directly get the multiplication matrix $S_k$ which must be a $d\times d$ integer matrix. Since we assumed $s_k^n=1,$ the inverse $s_k^{-1}$ is also an algebraic integer and thus $S_k^{-1}$ is an integer matrix. 

Independently of the basis, $P(\lambda)=0$ implies the matrix equation $P(L)=0.$ Since $P$ is minimal, $P$ is the characteristic polynomial of the matrix $L.$ Thus $\lambda$ is a complex eigenvalue of $L$ with a two-dimensional eigenspace, and we can always recover the complex number $z\in\QQ(\lambda )$ from the vector $\tilde{z}=x\in K$ by projection onto that eigenspace.

\subsection{The absence of numerical errors in the FT algorithm}\label{nonum} 
The IFS \eqref{alg} was lifted to an IFS on $K$ where $L$ and $S_k$ are the matrices representing multiplication by $\lambda$ and $s_k,$ respectively:
\begin{equation}
G(x)=Lx \quad \mbox{ and } \quad H_k(x)=S_k x + w_k \quad\mbox{ with } w_k=\tilde{v_k}\, .  
\label{lifted}\end{equation} 
We calculate on two levels. In $\CC$ we compute points of the attractor $A$ in order to draw a figure, or we determine the modulus of a complex number $z.$ These are numerical calculations with rounding errors.  The FT algorithm, however, is carried out in the ring of integers of $K$ which is isomorphic to $\ZZ^n .$ 
Actually, in all our examples, all coordinates fulfil $|x_j|<10.$  

The IFS maps $F_k=G^{-1}H_k$ on $K$ involve the inverse of the integer matrix $L.$
\begin{equation}  F_k(x)= L^{-1}(S_kx +w_k) \ .
\label{Fk}\end{equation}
Applying $F_k,$ we stay within rational numbers. When such calculations are iterated, however, the denominators become huge and we must resort to floating-point numbers. Thus for calculating iterated maps the space 
$K$ does not improve accuracy. For the FT algorithm, it does.

\begin{Proposition}\label{P6}
For an IFS \eqref{alg},\eqref{lifted} with algebraic integer data, the FT algorithm on $K$ works without numerical errors.
\end{Proposition}

{\it Proof. } Let $H(x)=Sx+w$ where $x,w$ are integer vectors in $K,$ and $S$ is a $d\times d$ integer matrix. The recursion step of the FT algorithm consists in determining $H'=F_j^{-1}HF_k$ for arbitrary $j,k\in\{ 1,...,m\}.$ With $F_j^{-1}(x)=S_j^{-1}(Lx-w_j)$ and \eqref{Fk} we get
\begin{equation}
F_j^{-1}HF_k(x)=S_j^{-1}\{ L(SL^{-1}(S_kx+w_k)+w)-w_j\}
=S_j^{-1}SS_kx +S_j^{-1}\{ Lw+Sw_k-w_j\}\, .
\label{FTK}\end{equation}
The matrices $S_k,S,L,S_j^{-1}$ all commute since they represent multiplication by complex numbers. So the calculation consists of vector additions and multiplications by integer matrices, and all results are integer vectors.

There is one further step in the FT algorithm: for every calculated neighbor map $H$ with translation vector $w=\tilde{v}$ we decide whether the modulus of the number $v$ exceeds the constant $C.$ If it does, the map $H$ cannot be proper. For this step we return from $K$ to the complex plane, by projection onto the $\lambda$-eigenspace of $L$ in the real vector space $\overline{K}.$ A $2\times d$ matrix for the projection $\pi$ is determined once, and for each new neighbor map $H(x)=Sx+w$  the modulus of $v=\pi(w)$ is calculated, with numerical errors.

However, the calculations are not iterated: we just check $|v|>C.$ And there is no harm if we incorrectly conclude $|v|\le C.$ This means that $H$ remains a candidate for a proper neighbor map until its successors lead to obviously greater values than $C,$ which happens exponentially fast (see the proof of Theorem \ref{equi}). So we need only care that rounding errors do not go upwards, or we can deliberately subtract a positive $\varepsilon$ from the calculated modulus. In this way, numerical errors never affect the result of the FT algorithm.
 \hfill $\Box$ \vspace{1ex}

\subsection{The case of Pisot factors}\label{mainproof}
Now we assume that the expansion factor $\lambda$ is a complex Pisot number. Since we assumed $s_k^n=1,$ the number of possible rotations $s$ in a neighbor map $h(z)=sz+v$ of our IFS is finite.  We prove Theorem \ref{main} by showing that  $K$ contains only finitely many integer vectors $w$ which can represent translations $v$ of proper neighbor maps.

The minimal polynomial $P(z)$ of $\lambda$ is the characteristic polynomial of the matrix $L.$ As a pair of Pisot numbers, $\lambda$ and $\overline{\lambda}$ have the largest modulus of all eigenvalues. The other eigenvalues, which we denote by $\lambda_3,...,\lambda_d,$ have modulus smaller than 1. We work in the real $d$-dimensional vector space $\overline{K}.$ Besides $B,$ which is a basis of $\overline{K},$ we consider a basis $E=\{ e_1,...,e_d\}$ of eigenvectors of $L.$ Let $e_1, e_2$ be basis vectors for the eigenplane of $\lambda ,$ and $e_3,...,e_d$ eigenvectors of the other eigenvalues.
  
The idea of proof is that the projection of the vector $w=\tilde{v}$ of a proper neighbor map $h(z)=sz+v$ to any eigenspace of $L$ must lie in a bounded set. Then the set of possible proper neighbor vectors $w$ is contained in a bounded set of $\overline{K}$ with respect to the basis $E,$ and also to basis $B.$ Boundedness in a finite-dimensional vector space does not depend on the basis. And a bounded set in  $\RR^n$ can contain only finitely many integer vectors.

For the eigenplane of $\lambda ,$ the boundedness of the projection $v$ of $w$ follows from the description of the FT algorithm in Section \ref{fin} and the calculation of $C$ in Section \ref{veri}.  We have $|v|\le C ,$ using the isomorphy between $\CC$ and the eigenplane.

Now we consider another eigenvalue $\lambda_q$ of $L,$ real or complex, with one- or two dimensional eigenspace $Y.$  The linear projection from $\overline{K}$ to $Y,$ which maps all other eigenspaces to zero, will be denoted by $\pi .$ The scalar product on $Y$ and the modulus function $|y|$ are determined up to a constant since $L$ acts on $Y$ as a similarity map with factor $|\lambda_q|.$  We fix the constant by requiring $|\pi (\tilde{1})|=1.$ Since $L$ and the symmetry matrices $S_j$ commute, $Y$ is an eigenspace of $S_j.$ In the Pisot case $S_j$ is a root of unity, so it acts as an isometry on $Y.$ 

We analyze the FT algorithm in Section \ref{fin}, as given in \eqref{FTK} on the integer points of $\overline{K},$ and the projections of the translation vectors $w$ to $Y.$ The algorithm starts with neighbor maps  $S_j^{-1}S_kx +S_j^{-1}(w_k-w_j)$ (this is \eqref{FTK} with $H=id$). In $Y$ we define
\[ c=2\cdot\max \{ |\pi (w_k)|\, |\, k=1,...,m\} .\]
Then
\[ |\pi (S_j^{-1}(w_k-w_j))|=  |\pi (w_k-w_j)| \le |\pi (w_k)|+|\pi (w_j)| \le c  \ .\]
Let us consider the recursion step. Some neighbor map $H(x)=Sx+w$  in \eqref{FTK} is given. The vector of a successor neighbor $H'$ has the form $w'=S_j^{-1}(Lw+Sw_k-w_j)\, .$  Thus
\[ |\pi(w')| = |\pi (Lw+Sw_k-w_j)| \le |\pi (Lw)|+|\pi(Sw_k)|+|\pi(w_j)| \le |\lambda_q|\cdot |w| + c \]
since $S$ is in the group generated by the $S_k$ acts as an isometry on $Y.$ We conclude that the moduli of all projections of neighbor vectors to $Y$ are bounded by the sum of a geometric series with factor $|\lambda_q|<1$ and starting value $c.$
\begin{equation} |\pi (w)|\le c/(1-|\lambda_q|)\ .
\label{LQ}\end{equation}
The set of projections of neighbor vectors to every eigenspace $Y$ of $L$ is bounded. This completes the proof of Theorem \ref{main}. Note that the estimate for the eigenspace of $\lambda ,$ given in \eqref{Ccon} with $r=1/|\lambda|,$ has the same form as \eqref{LQ}, with denominator $|\lambda|-1$ instead of $1-|\lambda_q|.$

\subsection{Comments on quantitative estimates}
Bandt and Mesing \cite{BM09} suggested to take the number of proper neighbor maps as a measure of complexity of the IFS, and to define a ``type'' as the relative position of two intersecting pieces up to similarity, given by a proper neighbor map.  Today we would say ``neighbor type''  since the number of neighborhoods, or the number of shapes of holes, or the maximum number of neighbors in a neighborhood, also indicate the complexity of the IFS.  Estimates for the number of proper neighbors of an IFS can tell us how much time the FT algorithm needs to check OSC and establish the neighbor graph.  The above proof does not provide estimates but gives an idea how different parameters influence the number of neighbor types. We combine this with the experience from our computer experiments.

The most important parameter is the degree $d$ of the Pisot number $\lambda .$ It is the dimension of $K,$ determines the number of eigenspaces $Y$ estimated by \eqref{LQ}, and thus influences the number of neighbor types in an exponential way.  We studied $d$ up to 8. For larger $d,$ the FT algorithm usually did not finish within an hour.

The group of possible rotations $s_k$ is the group of $n$-th root of unity. In the proof, $n$ is a factor of the number of neighbor types. What is worse, it is strongly connected with $d.$ If $n$ is odd (-1 not allowed as symmetry) then $d\ge n-1$ since $(z^n-1)/(z-1)$ is a minimal polynomial for the $s_k.$ If $n=2q$ is even, $d\ge q-1$ since we can take minimal polynomials $(z^q-1)/(z-1)$ and $(z^q+1)/(z+1).$  We studied odd $n$ up to 7 and even $n$ up to 14.

According to \eqref{LQ}, conjugates $\lambda_q$ with modulus near 1 will increase the number of possible neighbor vectors. Since we check the Pisot property of $\lambda$ by numerically determining the roots of the minimal polynomial, we can decide whether we like to work with some $\lambda$ with several conjugates near the unit circle.

The number $m$ of mappings in the IFS does not directly influence the number of neighbor vectors. However, the size of the $|\pi(w_k)|$ matters. Since the $v_k$ must be algebraic integers, we defined them directly as integer vectors $w_k=(x_1,...,x_d)$ in $K.$  We mostly took small coordinates $x_j$ between -3 and 3, for larger $d$ only between -1 and 1.  It is not clear how the choice of basis influences the number of neighbor types, and whether the size of $x_j$ has different effect for different $j.$

In most of our examples, the number of calculated neighbor vectors had the same magnitude as the estimated number of possible neighbor vectors in the proof of the theorem.  We had to severely restrict the degree of $\lambda ,$ the number of symmetries and the size of the $w_k$ in order to calculate the number of neighbor types and to draw the attractor.

\begin{figure}[h!t] 
\begin{center}
\includegraphics[width=0.38\textwidth]{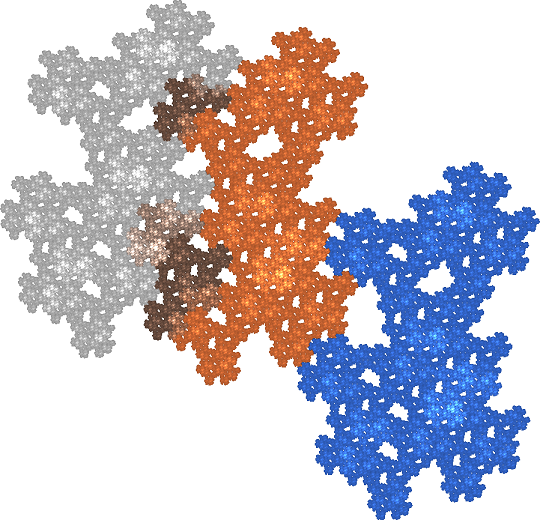} \ 
\includegraphics[width=0.58\textwidth]{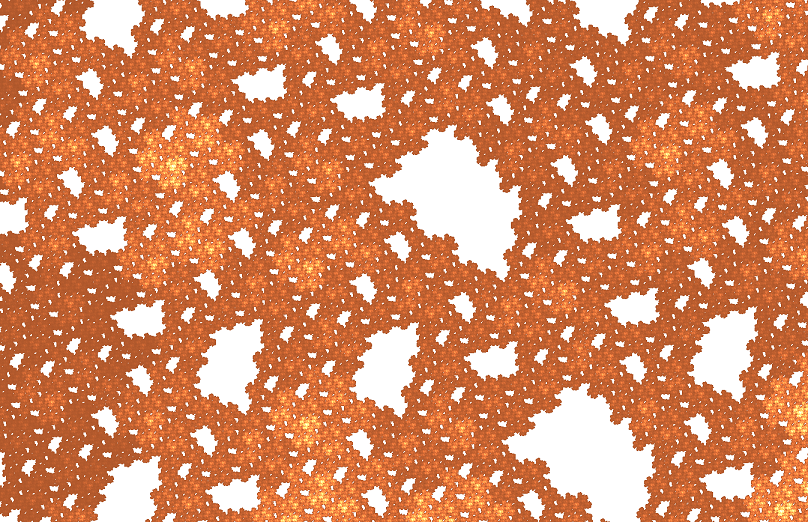} \vspace{1ex} \\
\includegraphics[width=0.38\textwidth]{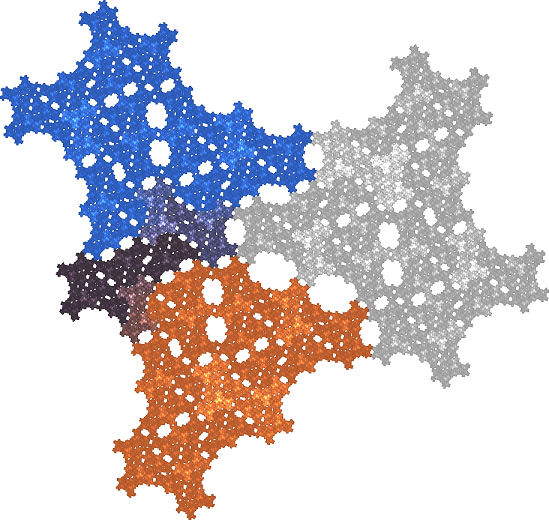} \
\includegraphics[width=0.58\textwidth]{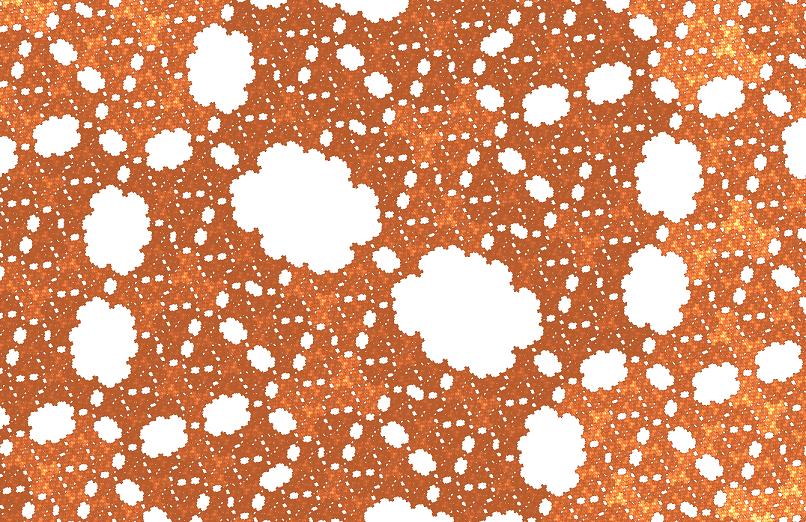} 
\end{center}
\caption{Two examples of low complexity from a Pisot number of degree 4.}\label{Fi5}
\end{figure}

\section{Examples} \label{details}

\subsection{Properties of our examples} 
\begin{figure}[h!t] 
\begin{center}
\includegraphics[width=0.29\textwidth]{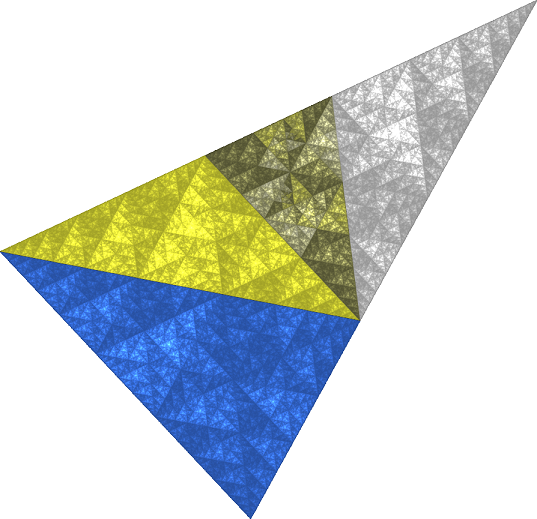} \
\includegraphics[width=0.26\textwidth]{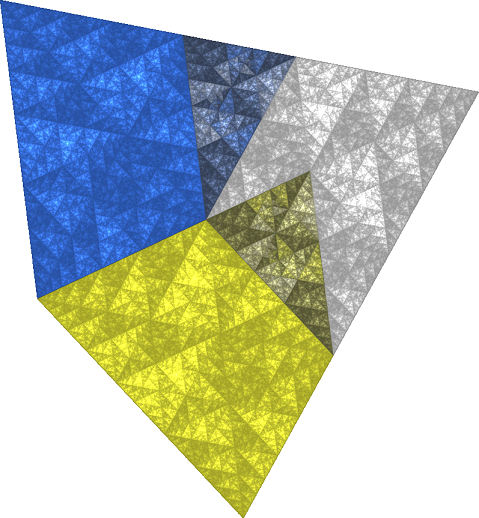} \
\includegraphics[width=0.36\textwidth]{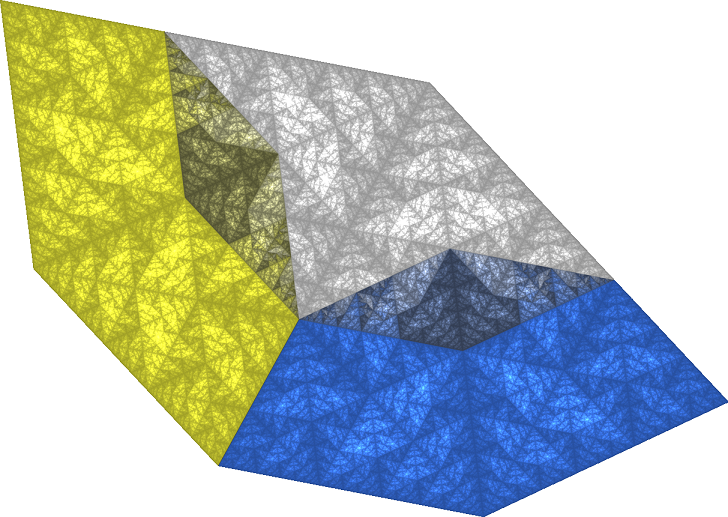} \vspace{1ex} \\
\includegraphics[width=0.4\textwidth]{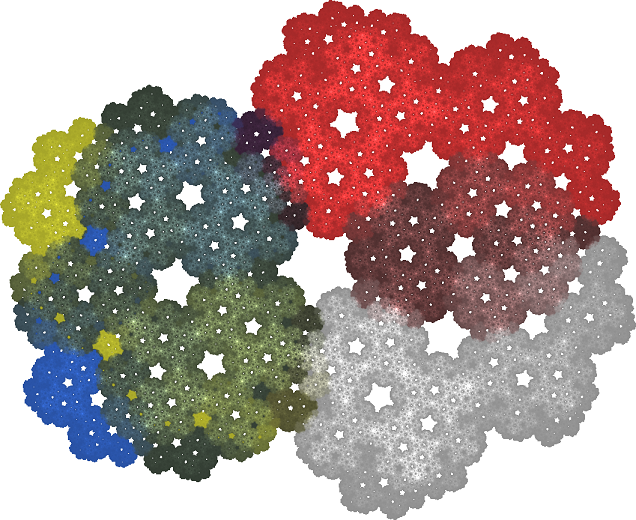} \
\includegraphics[width=0.55\textwidth]{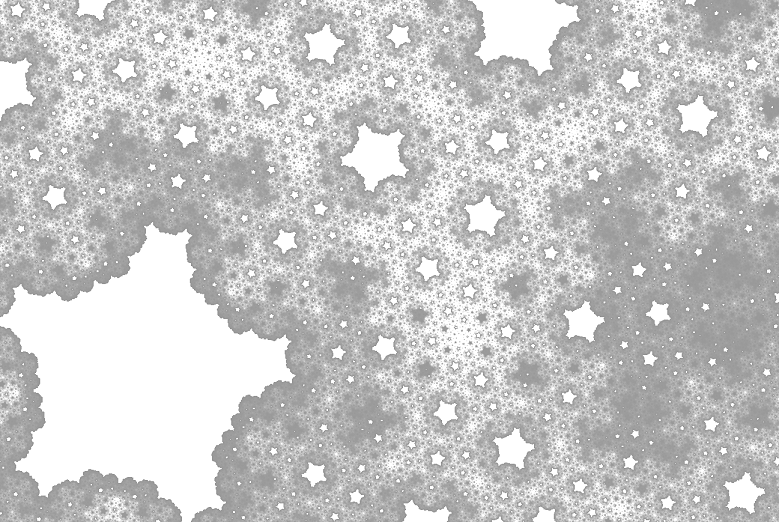} 
\end{center}
\caption{Examples from the Penrose family \cite{GS}. Upper row: The computer found overlapping schemes with three pieces for Robinson's triangle, Penrose's kite, and a double-kite. Lower row: Even with four pieces, there is a fractal pattern with very large overlaps.}\label{Fi7}
\end{figure}

For this chapter, we selected strongly connected attractors which represent carpets while trees or networks are known from papers on p.c.f. fractals \cite{Kig}. We avoided cutpoints, the removal of which disconnects a neighborhood of the point. They often lead to network-like structures. The dimension should be near to two, where the literature provides few examples. We got many overlapping tilings but mention only those connected with Penrose tiles in Figure \ref{Fi7}. Since we note the data of most IFS, we took small numbers of contractions. We required that the number of proper neighbor maps could be determined. 
Under these conditions, we present examples of very low complexity for further analysis of their neighbor graphs, and typical spaces with many different shapes of holes.  

The coloring of figures avoids artificial effects. The pieces in the global views are distinguished by different colors, but the local views are drawn with shades of one color. The shades indicate the overlaps which generate a structure of their own and play a part in the discussion of the pictures. Dark shades mean little or no overlap, very light shades stand for accumulated overlaps.

Let us briefly discuss degree $d=2.$
Any non-real quadratic integer $\lambda $ with modulus greater 1 is complex Pisot.
This holds in particular for Gaussian integers $\lambda = n+q\cdot i$ with symmetry group $\{ \pm 1, \pm i\} .$ OSC examples were studied in \cite{EFG1}. Another choice are the Eisenstein numbers $\lambda = n+q\cdot \omega$ with $\omega=\frac12 (1+\sqrt{3})$ and symmetry group $\{ \omega^k|\, k=1,...,6\} $ which were used for Figure \ref{Fi2}. Moreover, each quadratic integer $\lambda$ can be taken with trivial symmetries $\pm 1$ and integer translation vectors $v_k=a_k+b_k\lambda .$ In \cite{EFG2} such examples with irrational symmetries were studied.

\begin{figure}[h!t] 
\begin{center}
\qquad \includegraphics[width=0.278\textwidth]{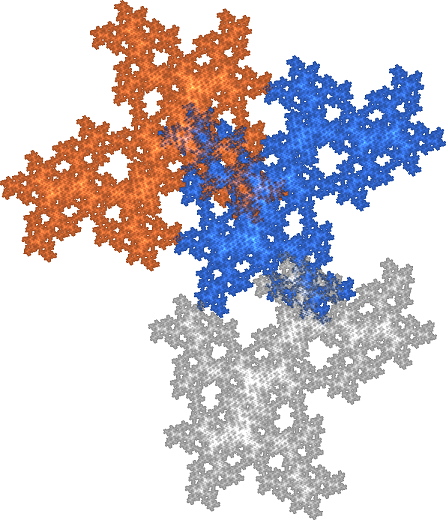} \quad\qquad\qquad\qquad
\includegraphics[width=0.49\textwidth]{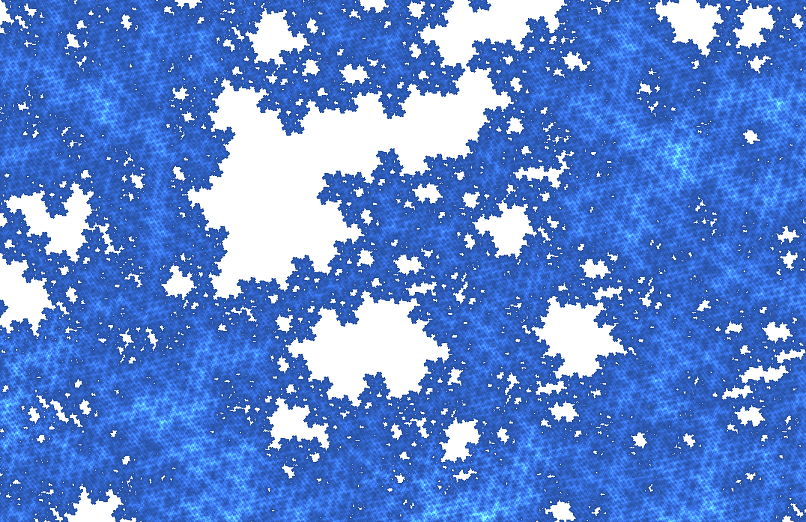} \vspace{.2ex} \\
\includegraphics[width=0.49\textwidth]{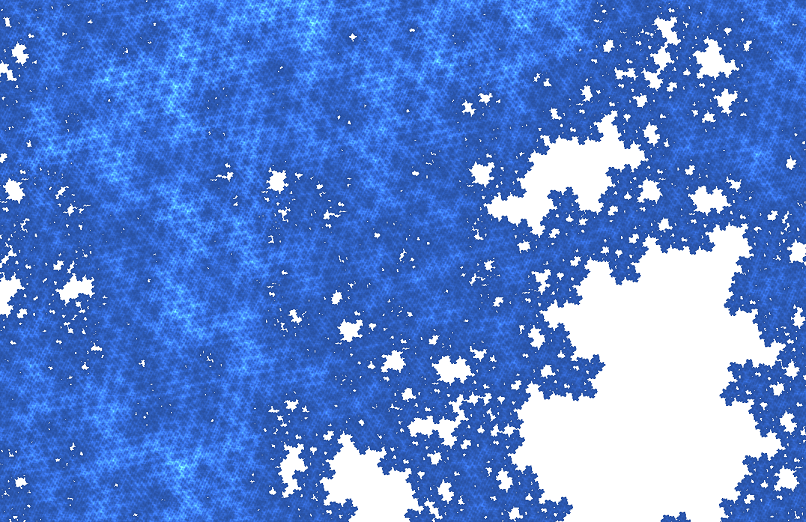} \
\includegraphics[width=0.49\textwidth]{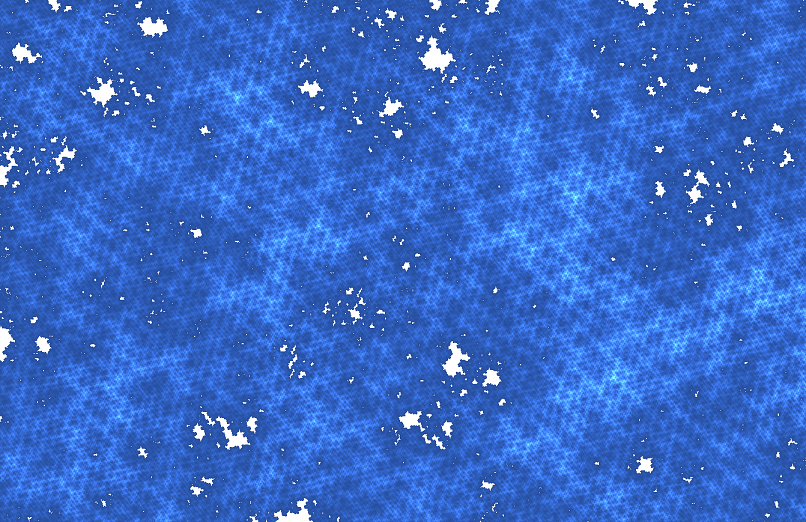} 
\end{center}
\caption{An example with 482 proper neighbor maps without rotation group from a complex Pisot number of degree 4.}\label{Fi6}
\end{figure}  

\subsection{Examples without rotation group} 
Suppose we want no more symmetries than $\pm 1.$ Then we only need a complex Pisot factor which can easily be found by a numeric search. A recent result of Bertin and Zaimi \cite{BZ} says that in any non-real algebraic number field $K,$ the complex Pisot numbers form a Delone set.   In the setting of Section \ref{sett}, we take the standard base $\{ 1,\lambda,...,\lambda^{d-1}\}$ and run Mekhontsev's program \cite{M} with random integer vectors $\tilde{v}=(x_1,...,x_d)$ to generate examples. 

We studied $d=4$ with palindromic polynomials $P(z)=z^4+az^3+bz^2+az+1.$ The inverse of a root of $P$ is a root, too. Thus if there are no real roots, we have a Pisot root. This will be true for $b>a^2/2$ since then $P(z)>z^2(z+\frac{a}{2})^2+(\frac{a}{2}z+1)^2\ge 0.$  

Our choice $a=2, b=4$ leads to $\lambda= -0.7429+1.5291\cdot i.$ The rotation angle of $\lambda$ is about $116^o$ and the modulus is $1.700.$ With $m=3$ maps, the formal similarity dimension is $\log 3/\log 1.7\approx 2.07.$ So there must be some overlap but it need not be large. In case of small overlaps, we should obtain examples with dimension near 2. 

In all computer experiments, we obtained many examples with high complexity and only few with a small number of neighbor maps. In Figure \ref{Fi5} there are two examples with small complexity, each with a single type of hole. In the upper row, all three pieces are translations of each other. All pieces of one level must be parallel to each other.  The pieces cannot be recognized, but for the holes this property is obvious. This example is an irregular analogue of the Sierpinski carpet. The data are $v_1=\lambda$ (grey), $v_2=\lambda+\lambda^2$ (orange), and $v_3=-1$ (blue).  There are 19 proper neighbor maps. 

The example in the bottom row of Figure \ref{Fi5} with a curious chain structure of holes was found by a search with larger $|v_k|.$ It includes the point reflection $s= -1$ but the holes are invariant under point reflection and again all parallel in the same level. The IFS is $h_1(z)=-(z+2+6\lambda+3\lambda^2+\lambda^3)$ (grey), 
$h_2(z)=z-2-\lambda$ (orange), and $h_3=z-\lambda+\lambda^2-\lambda^3$ (blue).  This IFS has only 24 proper neighbor maps, four of them with overlap. 

This example has a simple overlap structure. Two pairs of pieces on level 3 agree: $A_{211}=A_{333}$ and $A_{213}=A_{331}.$ The Hausdorff dimension $\beta$ can be determined by the method of Ngai and Wang. We give a much simpler argument. The $\beta$-dimensional Hausdorff measure is positive and finite, so let $\mu$ be the standardized Hausdorff measure. Then $\mu (A)=1$  and $\mu(A_k)=r^\beta$ with $r=1/|\lambda|$ for $k=1,2,3.$ Pieces of level 3 fulfil $\mu(A_w)=r^{3\beta}.$ Thus 
\[ 1= 3\cdot r^\beta - 2\cdot r^{3\beta}\]
since exactly two pieces are counted twice. Putting $x=r^\beta $ and dividing by $x-1$ we get the quadratic equation $2x^2+2x=1$ with solution $x=(\sqrt{3}-1)/2.$ The dimension is
$\beta =\log x /-\log |\lambda| \approx 1.894$ while the Sierpinski carpet satisfies $\alpha= \log 8 /\log 3 \approx 1.893.$

In both examples of Figure \ref{Fi5}, only two pieces overlap, and the holes are generated between the non-overlapping pairs of pieces which touch each other in a Cantor set. The overlaps just respect the hole structure for algebraic reasons. In Figure \ref{Fi6}, however, all pieces strongly overlap.  Various holes are generated between the blue piece and its neighboring pieces. On a higher level, these holes are modified by the overlaps in many ways. Only very few shapes of holes in the magnifications coincide. The IFS $h_1(z)=z+\lambda^2, \  h_2(z)=-(z-1+2\lambda^3)$ and $h_3(z)=z-\lambda$ does not look more complicated than those of Figure \ref{Fi5}. But the number of neighbor types is 482, and the structure is much more intricate. This was the typical situation in our experiments. 

\subsection{Examples from cyclotomic fields} 
For patterns with a symmetry group generated by a rotation $s$ of degree $d>2,$ it is difficult to represent $s$ in the standard basis of $\lambda .$ We better take the basis $\{ 1, s,..., s^{d-1}\}$ and search for Pisot numbers among the linear combinations of basis elements.  Figure \ref{Fi7} demonstrates the Penrose case with $s^5=-1$ (rotation by $36^o$ degrees) and $\lambda =s^3-1\approx -1.3090-0.9511 i.$  We first determined the companion matrix $S$ of the minimal polynomial $z^4-z^3+z^2-z+1$ of $s.$ Then $\lambda$ was obtained from the characteristic polynomial of $S^3-I.$ The Pisot property was checked numerically.  Actually, $\lambda= -\tau\cdot s$ where $\tau\approx 1.618$ is the golden mean.  

The celebrated aperiodic hat tile discovered recently \cite{hat} is a union of eight kites.  Our search with $m=3$ maps yielded three overlapping tiling schemes with the Robinson triangle, the Penrose kite \cite{GS} and a double-kite. The first two maps were the same for all three tiles: $h_1=s^4z, \  h_2=z+s^2 .$ The third is $h_2=sz+s^2, \ h_2=sz+s^3,$ and $h_2=s^{-2}z+1,$ respectively. The number of neighbor types is 86, 91 and 77, respectively. It is curious that overlapping schemes with one tile will be of finite type. It is not useful, however, since Penrose's kite-and-dart scheme and the two Robinson triangles \cite[chapter 11]{GS} are so much simpler.

The formal similarity dimension for 3 maps is $\log 3/\log\tau \approx 2.28.$ So we expect to obtain tiles for small overlaps. In contrast, the bottom row of Figure \ref{Fi7} was obtained with 4 maps and formal similarity dimension $\log 4/\log\tau \approx 2.88.$ Due to extremely large overlaps, the pattern is not plane-filling. The mappings are $h_1=s^3z, \ h_2=sz+s^3, \ h_3=s^7(z-1),$ and $h_4=-z.$ Here we had 55 proper neighbors.

\begin{figure}[h!t] 
\begin{center}
\ \includegraphics[width=0.485\textwidth]{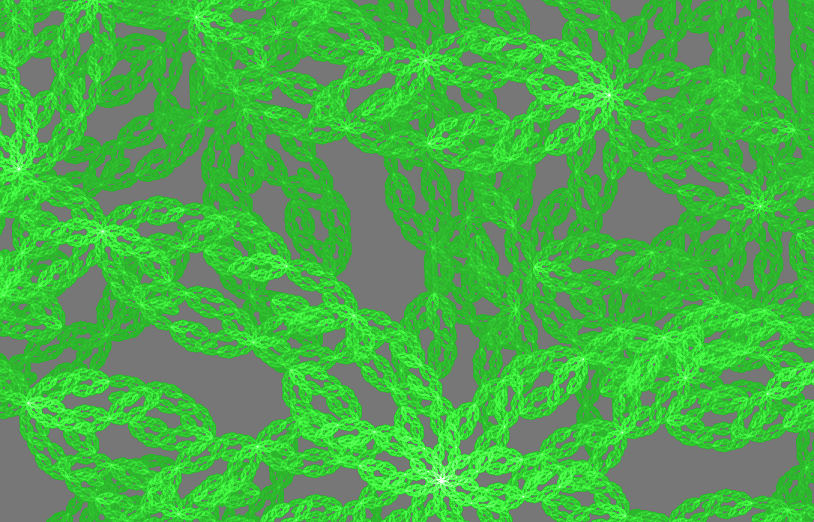} \
\includegraphics[width=0.485\textwidth]{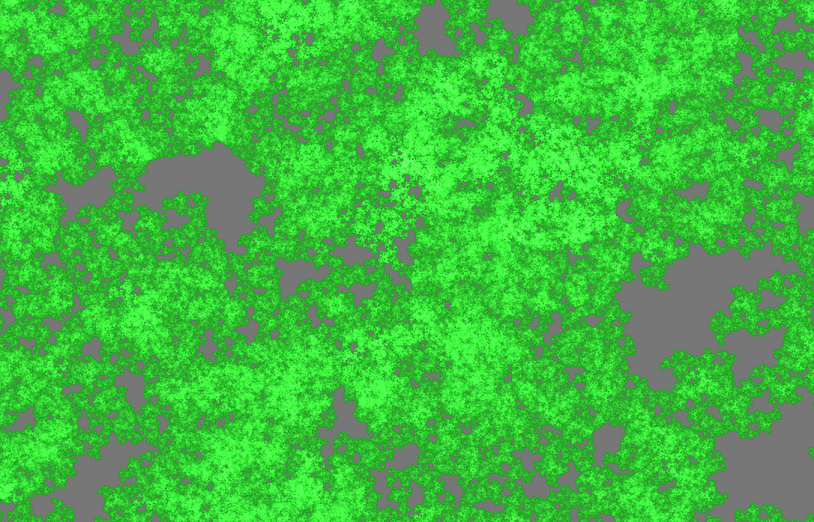} \vspace{1ex} \\
\includegraphics[width=0.485\textwidth]{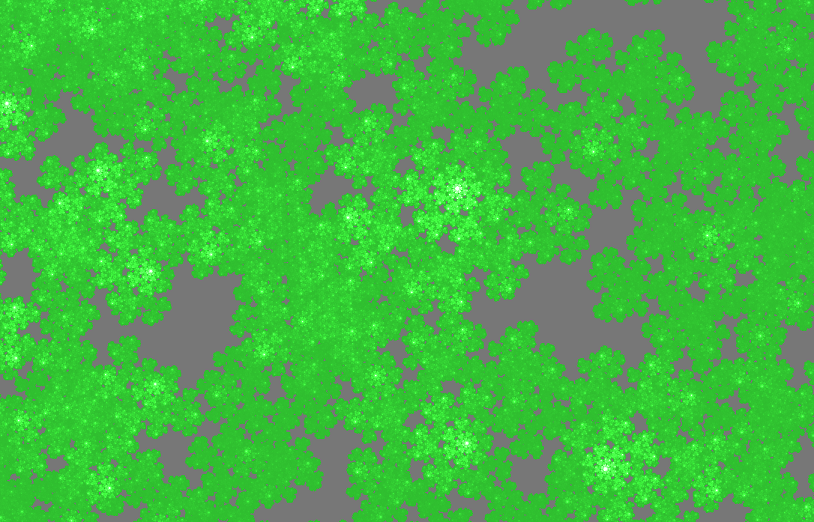} \
\includegraphics[width=0.485\textwidth]{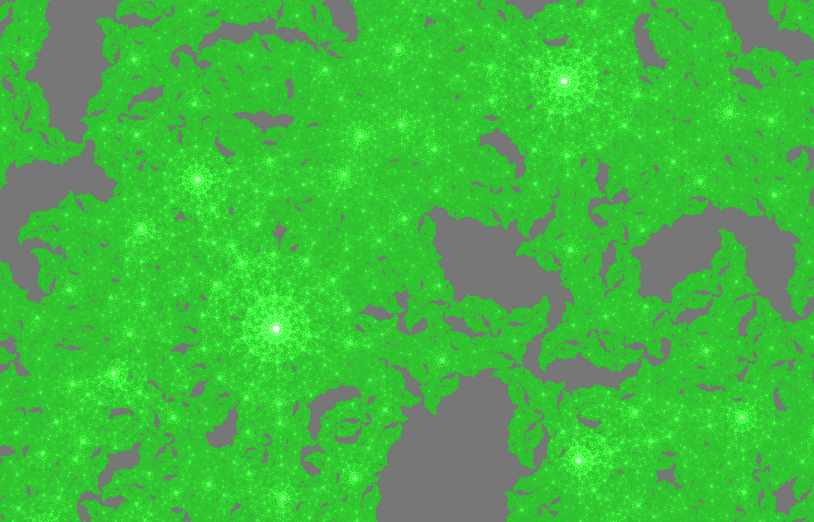} 
\end{center}
\caption{Local views of four IFS with sevenfold symmetry: $s^7=-1$ and $\lambda=s^{-1}+s+s^3.$ Upper row with $m=5$ and small overlap, 703 and 2284 neighbor types. Lower row with $m=6$ maps and larger overlap, 989 and 1776 neighbor types. }\label{Fi8}
\end{figure}

Another Pisot number is $\lambda= 1-s^3-s^4=\tau^2s.$ With $m=7$ maps the formal similarity dimension is $\log 7/\log\tau^2 \approx 2.02.$ This leads to Figure \ref{Fi1} with 1977 neighbors and almost plane-filling parts. 
The examples in Figure \ref{Fi8} with $m=6$ have smaller dimension. 

Sevenfold symmetry was already studied in Figure \ref{Fi3}, with $s^7=-1$ and a non-Pisot factor. For Figure \ref{Fi8} we take $\lambda=s^{-1}+s+s^3\approx 2.0245-0.9749i$ which is a Pisot unit, and show one local view for four different IFS of high complexity. All have their specific geometry of holes and appearance of symmetry.
 
\begin{figure}[h!t] 
\begin{center}
\qquad\ \includegraphics[width=0.3\textwidth]{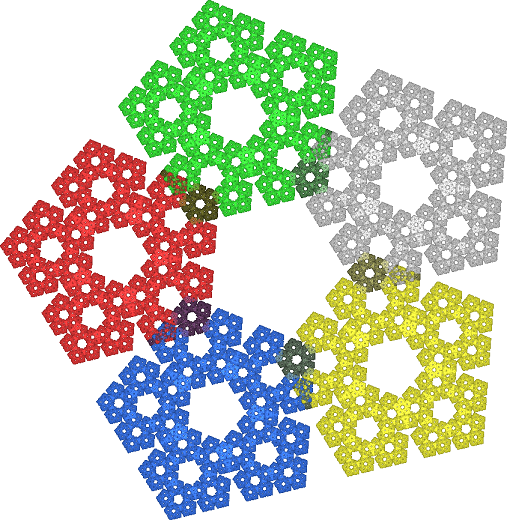} \quad\qquad\qquad\
\includegraphics[width=0.485\textwidth]{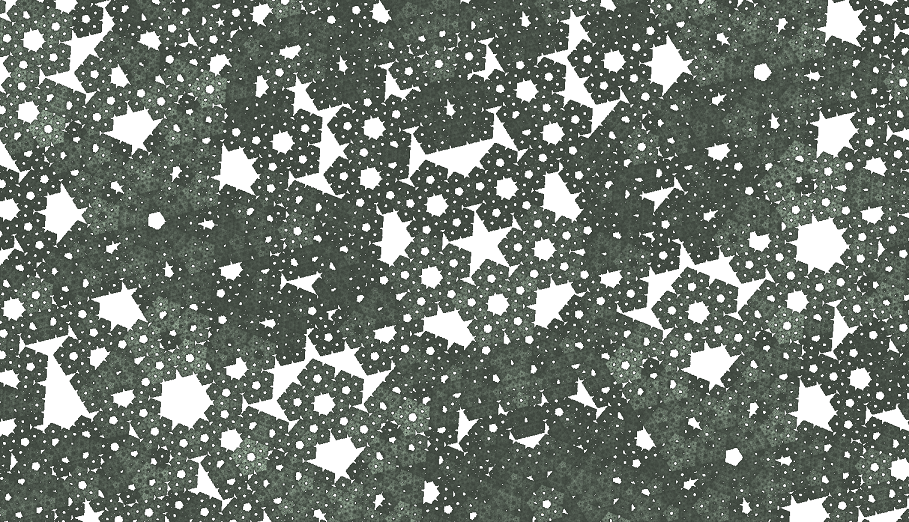} \vspace{1ex} \\
\includegraphics[width=0.485\textwidth]{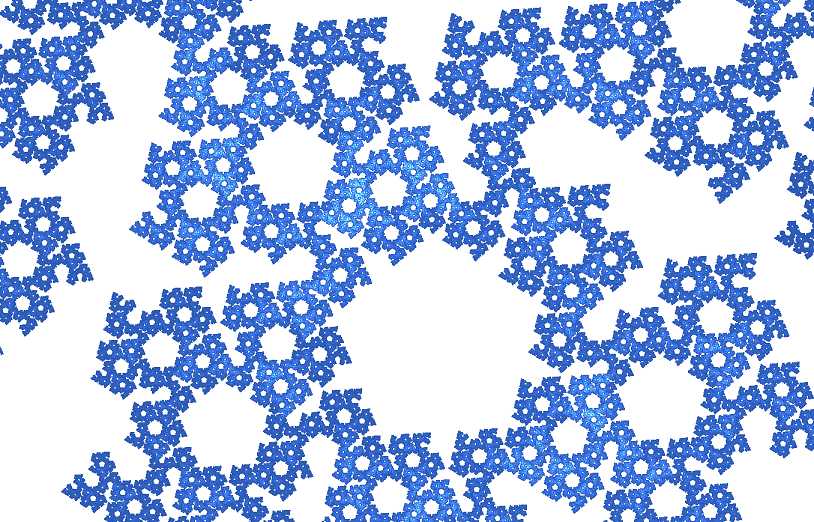} \
\includegraphics[width=0.485\textwidth]{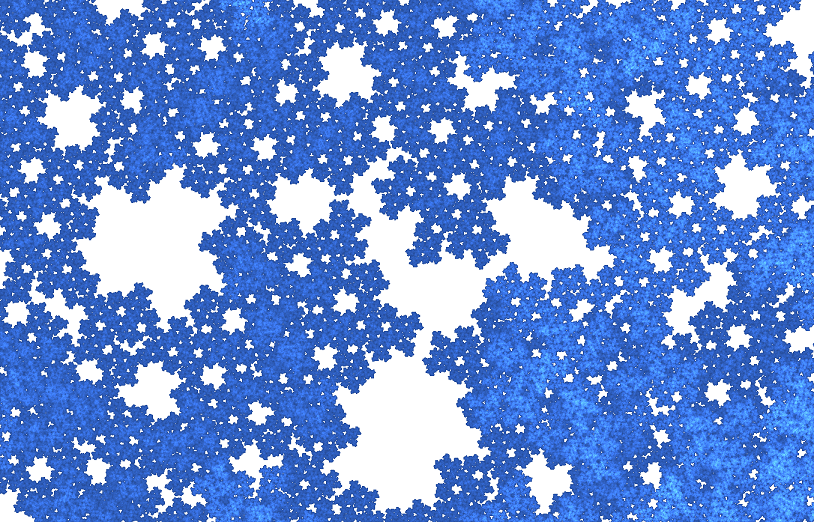} \vspace{1ex} \\
\includegraphics[width=0.485\textwidth]{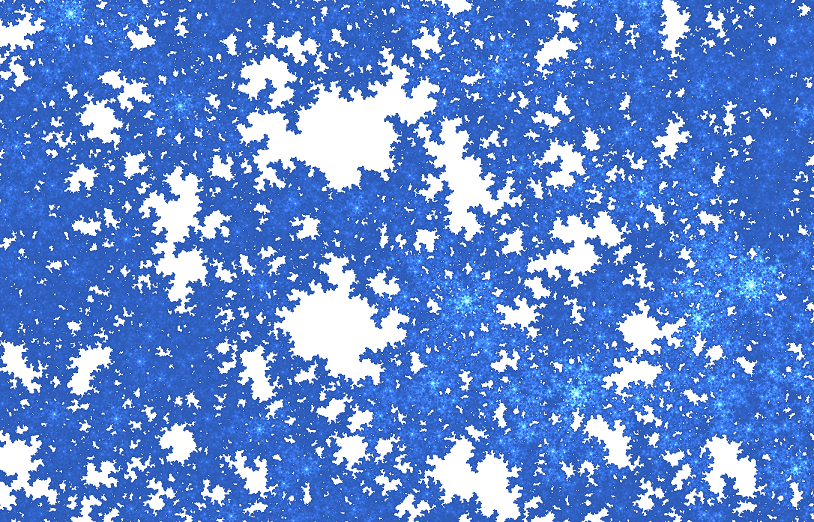} \
\includegraphics[width=0.485\textwidth]{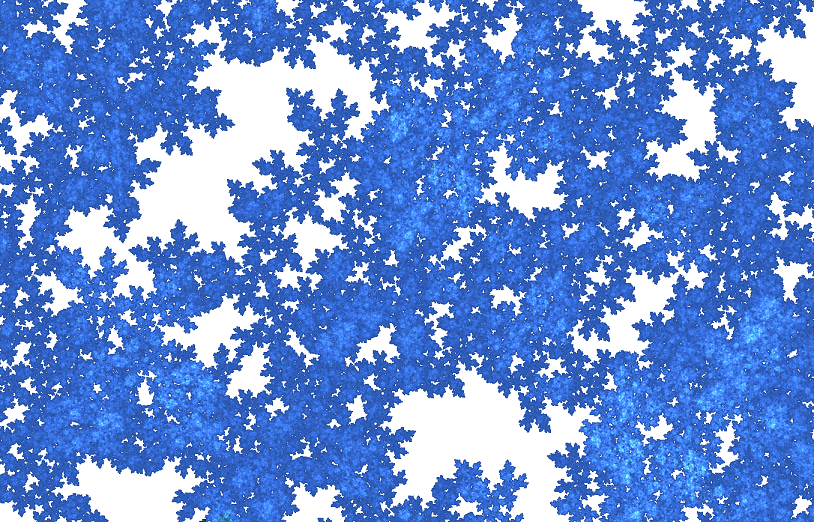} 
\end{center}
\caption{Top row: A complex Pisot factor was determined so that certain pieces of the fractal pentagon coincide. Many other partial overlaps and new shapes of holes are generated, raising the number of neighbor types to 895. Below: A search with this Pisot factor gives few simple and many complicated examples. We show local views for an IFS with 13 and three IFS with about 6000 neighbor types.}\label{Fi9}
\end{figure}  

\subsection{Examples from prescribed overlaps} \label{prescribed}
IFS for patterns with symmetries can also be determined by prescribing exact overlaps. This is demonstrated for an IFS with $h_k=s^kz$ for $k=0,...,d-1$ with full symmetry of degree $d.$ We calculate $\lambda$ so that two pieces of third level coincide. For degree 3 and 4, this leads to the golden gasket \cite{BMS} and to the overlapping square in \cite{BM09}. In Figure \ref{Fi9}, we take $d=5$ and determine a real $\lambda >1 .$ Let $r=1/\lambda$ and $s^5=1.$ The coincidence $A_{022}= A_{144}$ is expressed by the equation
 \[ 1 +rs^2 +r^2s^2 - (s +rs^4 +r^2s^4)=0 \ .\]
Dividing by $1-s$ and multiplying by $\lambda^2$ we get
\[ \lambda^2 + (s^2+s^3)(\lambda +1) =0 \ .\] 
This is a quadratic equation for $\lambda$ with coefficients in the ring generated by $s.$ Thus we have a quadratic extension of the 5th cyclotomic field, and $\lambda$ will be an algebraic integer of degree $2\cdot 4=8.$ It is convenient to take the basis $\{ 1,s,s^2,s^3, \lambda, \lambda s, \lambda s^2, \lambda s^3\}$ for which the multiplication with $\lambda$ can be directly expressed by an integer matrix $L.$ The characteristic polynomial of $L$ is  
\[ P(z)=z^8-2z^7+7z^6-5z^5+6z^4-5z^3+3z^2-z+1 \]
with double root  $\lambda \approx 2.3165$ and $|\lambda_q|\le 0.79$ for the conjugates. We replace $\lambda$ by $s\lambda$ to get a complex Pisot number.  As Figure \ref{Fi9} indicates, the overlap generates other partial overlaps and lots of different holes in the magnifications. However, Theorem \ref{main} says that the number of neighbor types, and also of shapes of holes, is finite. The FT algorithm calculated 895 neighbor types. A search with the Pisot expansion factor and randomly modified $h_k$ yields few patterns with low complexity and many with several thousand neighbor types, as illustrated in  Figure \ref{Fi9}.
The ``cactus carpet'' clearly has very small dimension and complexity. For the other three attractors we chose views which look similar in dimension and density of holes. Still, we can distinguish their geometry by eyesight. Can we express this with mathematical parameters?  

\subsection{Outlook} 
A wide variety of self-similar plane sets with positive Hausdorff measure in fractal dimensions near to 2 was uncovered. Overlaps form a mechanism to generate many natural-looking fractally homogeneous subsets of the plane.  So far, the overlaps are vaguely controlled by the data. For complex examples, algorithms for drawing local views can run out of time. It may be impossible to count neighbors and neighbourhoods and to calculate Hausdorff dimension. 

New challenges come up. How can we describe the geometry in more complex examples? Are overlaps necessary to build such sets? Can we generalize the present results to the non-commutative setting in dimension three? Can we construct textures with given properties? 
Fractal models shall contribute to understanding nature, as envisioned by Mandelbrot \cite{FGN} and Barnsley \cite{Bar}. From satellite images up to microscopic views in medicine, data now show more fractal features than 30 years ago.

\end{document}